\documentclass[12pt,a4paper]{amsart}
\usepackage{amsmath,amsthm,amssymb,amscd}
\usepackage[arrow,matrix]{xy}
\usepackage{mathtools}
\usepackage{enumerate}
\usepackage{ascmac}
\usepackage{graphicx}
\usepackage{mathrsfs}
\usepackage{typearea}
\usepackage{ytableau}
\usepackage{braket}
\usepackage{caption}
\usepackage{tikz}
\usetikzlibrary{calc}
\usepackage{xcolor}
\usepackage[draft]{changes}
\usetikzlibrary{decorations.pathmorphing}
\usepackage{thmtools}
\usepackage{thm-restate}

\calclayout

\theoremstyle{definition}
\newtheorem{thm}{Theorem}[section]
\newtheorem{conj}[thm]{Conjecture}
\newtheorem{lemma}[thm]{Lemma}
\newtheorem{corollary}[thm]{Corollary}
\newtheorem{remark}[thm]{Remark}
\newtheorem{dfn}[thm]{Definition}

\newtheorem{prop}[thm]{Proposition}

\def\Span{\mathop{\mathrm{Span}}\nolimits}

\def\dim{\mathop{\mathrm{dim}}\nolimits}

\def\Ker{\mathop{\mathrm{Ker}}\nolimits}

\def\det{\mathop{\mathrm{det}}\nolimits}

\def\GL{\mathop{\mathit{GL}}\nolimits}

\newcommand{\mf}[1]{{\mathfrak{#1}}}
\newcommand{\mb}[1]{{\mathbf{#1}}}
\newcommand{\bb}[1]{{\mathbb{#1}}}
\newcommand{\mca}[1]{{\mathcal{#1}}}
\newcommand{\mr}[1]{{\mathrm{#1}}}
\newcommand{\msc}[1]{{\mathscr{#1}}}
\newcommand{\msf}[1]{{\mathsf{#1}}}
\newcommand{\relmiddle}[1]{\mathrel{}\middle#1\mathrel{}}
\newcommand{\da}[1]{\big\downarrow\raise.5ex\rlap{$\scriptstyle#1$}}
\newcommand{\Prob}{\mathbb{P}}

\DeclareMathOperator{\Res}{\msf{Res}}

\makeatletter
\@namedef{subjclassname@2020}{\textup{2020} Mathematics Subject Classification}
\makeatother

\begin{document}

\title{A proof of the Stanley--Stembridge conjecture}
\author{Tatsuyuki Hikita}
\address{\textsc{Research Institute for Mathematical Sciences, Oiwake Kita-Shirakawa Sakyo Kyoto 606-8502 JAPAN}}
\email{\texttt{thikita@kurims.kyoto-u.ac.jp}}
\date{}

\begin{abstract}
	We give a probabilistic interpretation of the coefficients of the elementary symmetric function expansion of the chromatic quasisymmetric function for any unit interval graph. As a corollary, we prove the Stanley--Stembridge conjecture.
\end{abstract}

\maketitle

\section{Introduction}

The Stanley--Stembridge conjecture is a long-standing conjecture in algebraic combinatorics stated by Stanley and Stembridge \cite{SS93} in 1993, and reformulated by Richard Stanley \cite{Sta95} in terms of the notion of chromatic symmetric functions in 1995. This conjecture is also known as the 21st problem in Stanley's list of positivity problems \cite{Sta00}. In this paper, we prove this conjecture affirmatively by exhibiting a probabilistic formula for the coefficients of the elementary symmetric function expansion ($e$-\emph{expansion} for short) of the chromatic quasisymmetric function for any unit interval graph.

\subsection{Main results}\label{Sec:Main}

Let $\Gamma$ be a unit interval graph with $n$ vertices (see Section~\ref{Sec:UIG} for more details). We write the $e$-expansion of its chromatic quasisymmetric function $\mb{X}_{\Gamma}(x;q)$ (see Definition~\ref{Dfn_CQF}) as
\begin{align}\label{Eqn:eCoeffDef}
	\mb{X}_{\Gamma}(x;q)=\sum_{\lambda\vdash n}c_{\lambda}(\Gamma;q)e_{\lambda}(x).
\end{align}
Our starting point is the following observation\footnote{When $q=1$, this formula easily follows from the definition of the chromatic symmetric function by comparing the coefficients of $x_1x_2\cdots x_n$ in (\ref{Eqn:eCoeffDef}). Indeed, the coefficient of $x_1x_2\cdots x_n$ in $\mb{X}_{\Gamma}(x;1)$ is $n!$ and the coefficient in $e_{\lambda}(x)$ is $\frac{n!}{\prod_i\lambda_i!}$. We note that this argument works for any simple graph.} (see Corollary~\ref{Cor:Prob_Rel}):
\begin{align*}
	\sum_{\lambda\vdash n}\frac{c_{\lambda}(\Gamma;q)}{\prod_{i}[\lambda_i]_q!}=1.
\end{align*}

In light of the refined Stanley--Stembridge conjecture (see Conjecture~\ref{Conj_SW}) due to Shareshian and Wachs, this equality suggests that the map $\msf{P}_{\Gamma}$ from the set of all partitions $\msf{Par}$ to $\bb{R}$ sending $\lambda\in\msf{Par}$ to
\begin{align*}
	\msf{P}_{\Gamma}(\lambda;q)\coloneqq\begin{cases}
		\frac{c_{\lambda}(\Gamma;q)}{\prod_{i}[\lambda_i]_q!} & \text{if }\lambda\vdash n,\\
		0 & \text{otherwise,}
	\end{cases}
\end{align*}
should define a probability measure on $\msf{Par}$, i.e., $\msf{P}_{\Gamma}(\lambda;q)\geq0$ for any $\lambda$ and $q\in\bb{R}_{>0}$. Our main result is to prove this inequality in general, which is sufficient to conclude that the original Stanley--Stembridge conjecture is true. 

\subsubsection{Conjugate Hessenberg functions}

For each $i\in[n]$, let $\Gamma_{\leq i}$ be the subgraph of $\Gamma$ obtained by removing the vertices $\{i+1,\ldots,n\}$ and edges incident to them. We note that the unit interval graph $\Gamma_{\leq i}$ is determined by its subgraph $\Gamma_{\leq i-1}$ and the value
\begin{align*}
	\msf{e}(i)\coloneqq\min\{j\in\{1,\ldots,i\}]\mid \{j,i\}\in\Gamma_{\leq i}\text{ or }j=i\}-1,
\end{align*}
since $\{j,i\}$ is an edge of $\Gamma_{\leq i}$ if and only if $\msf{e}(i)<j<i$. In particular, $\Gamma$ is determined by the sequence $\msf{e}=(\msf{e}(1),\ldots,\msf{e}(n))$ which we call the \emph{conjugate Hessenberg function} of $\Gamma$. Conversely, a sequence $\msf{e}=(\msf{e}(1),\ldots,\msf{e}(n))$ comes from a unit interval graph if and only if
\begin{align*}
	0\leq\msf{e}(i)<i,\quad\msf{e}(i)\leq\msf{e}(i+1)
\end{align*}
for any $i\in\{1,\ldots,n\}$, with the understanding that $\msf{e}(n+1)=\infty$ (see Section~\ref{Sec:conj_Hess}).

\subsubsection{Probabilistic ansatz}

We prove the positivity of $\msf{P}_{\Gamma}(\lambda;q)$ by constructing an explicit stochastic process whose probability distribution realizes $\msf{P}_{\Gamma}$. We first explain how to find such a stochastic process by borrowing probabilistic terminology. We refer readers who are not familiar with probability theory to Appendix~\ref{Sec:Prob} or any standard textbooks for the relevant terminology. We also refer to \cite{Hik25} for details about the actual experiments we performed.

For any sequence of integers $\msf{e}:\bb{Z}_{>0}\rightarrow\bb{Z}$ (not necessarily a conjugate Hessenberg function), we want to construct a stochastic process $\{\msf{X}_n^{(\msf{e})}\}_{n\in\bb{Z}_{\geq0}}$ with state space $\msf{Par}$, depending on $\msf{e}$ and $q\in\bb{R}_{>0}$, such that if $\msf{e}_{\leq n}\coloneqq(\msf{e}(1),\ldots,\msf{e}(n))$ is a conjugate Hessenberg function corresponding to a unit interval graph $\Gamma_{\leq n}$, then 
\begin{align*}
	\Prob\{\msf{X}_n^{(\msf{e})}=\lambda\}=\msf{P}_{\Gamma_{\leq n}}(\lambda;q).
\end{align*}

We look for a collection of maps $\msf{K}_{n}^{(\msf{e})}:\msf{Par}^{n+1}\rightarrow[0,1]$ satisfying 
\begin{align*}
	\Prob\{\msf{X}^{(\msf{e})}_{0}=\lambda_0,\ldots,\msf{X}^{(\msf{e})}_{n}=\lambda_n\}=\msf{K}^{(\msf{e})}_0(\lambda_0)\msf{K}^{(\msf{e})}_1(\lambda_0;\lambda_1)\cdots\msf{K}^{(\msf{e})}_n(\lambda_0,\ldots,\lambda_{n-1};\lambda_n),
\end{align*}
by imposing the following simplifying assumptions:
\begin{enumerate}
	\item The stochastic process $\{\msf{X}_n^{(\msf{e})}\}_{n\in\bb{Z}_{\geq0}}$ is a growth process of Young diagrams, i.e., $\Prob\{\msf{X}^{(\msf{e})}_{0}=\lambda_0,\ldots,\msf{X}^{(\msf{e})}_{n}=\lambda_n\}\neq0$ only when $|\lambda_i|=i$ and $\lambda_0\subset\lambda_1\subset\cdots\subset\lambda_n$.
	\item The transition rule from time $n-1$ to time $n$ depends only on the number $\msf{e}(n)$, i.e., there exists a map $\msf{K}^{(r)}_{n}:\msf{Par}^{n+1}\rightarrow[0,1]$ for each $r\in\bb{Z}_{\geq0}$ and $n\in\bb{Z}_{>0}$ such that $\msf{K}^{(\msf{e})}_{n}=\msf{K}^{(\msf{e}(n))}_{n}$ for any integer sequence $\msf{e}$.
\end{enumerate}
In particular, the random variable $\msf{X}^{(\msf{e})}_{0}$ must be the constant map to the empty partition $\emptyset$.

We note that the data $(\lambda_0,\ldots,\lambda_n)$ satisfying $|\lambda_i|=i$ and $\lambda_0\subset\lambda_1\subset\cdots\subset\lambda_n$ is equivalent a standard Young tableau of size $n$. Hence we may regard $\msf{K}^{(r)}_{n}$ as a transition probability $\msf{P}^{(r)}_{n}:\msf{SYT}_{n-1}\times\msf{SYT}_{n}\rightarrow[0,1]$ from the set $\msf{SYT}_{n-1}$ of standard Young tableaux of size $n-1$ to $\msf{SYT}_{n}$, defined by
\begin{align*}
	\msf{P}^{(r)}_n(T;T')=\begin{cases}
		\msf{K}^{(r)}_n(\lambda_0,\ldots,\lambda_{n-1};\lambda_n) & \text{if }T\subset T',\\
		0 & \text{otherwise},
	\end{cases}
\end{align*}
where $T$ (resp. $T'$) corresponds to $\lambda_0\subset\cdots\subset\lambda_{n-1}$ (resp. $\lambda_0\subset\cdots\subset\lambda_{n}$). In particular, the stochastic process $\{\msf{X}_n^{(\msf{e})}\}_{n\in\bb{Z}_{\geq0}}$ can be lifted to a time-inhomogeneous Markov chain $\{\widetilde{\msf{X}}_n^{(\msf{e})}\}_{n\in\bb{Z}_{\geq0}}$ with state space $\msf{SYT}\coloneqq\cup_n\msf{SYT}_n$ and transition probabilities $\{\msf{P}^{(\msf{e}(n))}_{n}\}_{n\in\bb{Z}_{>0}}$.

Surprisingly, from these assumptions and the observable data  $\{\msf{P}_{\Gamma}(\lambda;q)\}_{\lambda,\Gamma}$, one can determine the data $\{\msf{P}^{(r)}_{n}\}$ step by step in a consistent manner, i.e., without violating the simplifying assumptions. Through numerical experiments, we were able to guess a formula for $\msf{P}^{(r)}_{n}$. 

\subsubsection{Maya diagrams}

In order to describe our formula for $\msf{P}^{(r)}_{n}$, we first construct a map 
\begin{align*}
	\boldsymbol{\delta}^{(r)}:\msf{SYT}\rightarrow\msf{Maya}
\end{align*}
from $\msf{SYT}$ to the set $\msf{Maya}$ of \emph{Maya diagrams}, i.e., infinite sequences $\delta=(\delta_i)_{i\in\bb{Z}}$ of $0$ or $1$ such that $\delta_i=0$ for $i\gg0$ and $\delta_i=1$ for $i\ll0$. We will depict a Maya diagram by a sequence of white and red boxes, where $\delta_i=0$ (resp. $\delta_i=1$) corresponds to the white (resp. red) box.

\begin{dfn}\label{Def:Delta}
Let $T\in\msf{SYT}$ be a standard Young tableau and $r\in\bb{Z}_{\geq0}$. We define the Maya diagram $\boldsymbol{\delta}^{(r)}(T)=(\delta_i)_{i\in\bb{Z}}\in\msf{Maya}$ as follows:
\begin{itemize}
	\item If $i\leq0$, then we set $\delta_i=1$.
	\item If $i$ is greater than the length of the first row of $T$, then we set $\delta_i=0$.
	\item Otherwise, let $\msf{top}_T(i)$ be the entry in the topmost box of the $i$-th column of $T$ (in the French convention). If $\msf{top}_T(i)>r$, then we set $\delta_i=1$. Otherwise, we set $\delta_i=0$.
\end{itemize}
\end{dfn}

In other words, $\boldsymbol{\delta}^{(r)}(T)$ is obtained by first looking at $T$ from above, coloring the box $\boxed{i}$ with $i>r$ red, and then adding infinitely many red boxes to the left and infinitely many white boxes to the right. 

\begin{dfn}\label{Def:varphi}
Let $\delta=(\delta_i)_{i\in\bb{Z}}\in\msf{Maya}$ be a Maya diagram. We set
\begin{align*}
	W(\delta)&\coloneqq\{i\in\bb{Z}\mid \delta_i=0,\delta_{i-1}=1\},\\
	R(\delta)&\coloneqq\{i\in\bb{Z}\mid \delta_i=1,\delta_{i-1}=0\}.
\end{align*}
For any $c\in \bb{Z}$, we define $\varphi_c(\delta;q)\in\bb{Q}(q)$ (or simply $\varphi_c(\delta)$) by
\begin{align*}
	\varphi_c(\delta;q)\coloneqq\begin{cases}
		\frac{\prod_{b\in R(\delta)}[c-b]_q}{\prod_{a\in W(\delta)\setminus\{c\}}[c-a]_q} &\text{if }c\in W(\delta),\\
		0 &\text{otherwise.}
	\end{cases}
\end{align*}
\end{dfn}

In other words, $W(\delta)$ (resp. $R(\delta)$) is the set of indices of the leftmost white (resp. red) boxes in each run of consecutive white (resp. red) boxes. See Figure~\ref{Fig:Delta} for an example.

\begin{figure}[h]
\begin{tikzpicture}
\def\sz{1.5em}
\node at (0,0) {$T=$}; 
\node at (3*\sz,0) {\begin{ytableau}[] 
4 & *(red) 5 \\
1  & 2 & 3  & *(red) 6
\end{ytableau}};
\filldraw[fill=red] (1*\sz, -4*\sz) -- (1*\sz, -3*\sz) -- (0.5*\sz, -3*\sz)
    decorate [decoration={snake, amplitude=1pt, segment length=4pt}] { -- (0.5*\sz, -4*\sz) }
    -- cycle;
\node at (-0.25*\sz, -3.5*\sz) {$\cdots$};
\node at (-1.5*\sz,-3.5*\sz) {$\delta=$};
\filldraw[fill=white] (5*\sz, -4*\sz) -- (5*\sz, -3*\sz) -- (5.5*\sz, -3*\sz)
    decorate [decoration={snake, amplitude=1pt, segment length=4pt}] { -- (5.5*\sz, -4*\sz) }
    -- cycle;
\node at (6.25*\sz, -3.5*\sz) {$\cdots$};
\draw[->,thick] (3*\sz,-1.5*\sz) --node[right] {$\boldsymbol{\delta}^{(4)}$} (3*\sz, -2.5*\sz);
\draw[fill=white] (1*\sz, -4*\sz) rectangle ++(\sz, \sz);
\draw[fill=red]   (2*\sz, -4*\sz) rectangle ++(\sz, \sz);
\draw[fill=white] (3*\sz, -4*\sz) rectangle ++(\sz, \sz);
\draw[fill=red]   (4*\sz, -4*\sz) rectangle ++(\sz, \sz);
\node[font=\Large] at (11.25*\sz,0.25*\sz) {$\varphi_{1}(\delta)=\frac{[-1]_q[-3]_q}{[-2]_q[-4]_q}$};
\node[font=\Large] at (11*\sz,-1.75*\sz) {$\varphi_{3}(\delta)=\frac{[1]_q[-1]_q}{[2]_q[-2]_q}$};
\node[font=\Large] at (10.75*\sz,-3.75*\sz) {$\varphi_{5}(\delta)=\frac{[3]_q[1]_q}{[4]_q[2]_q}$};
\end{tikzpicture}
\caption{Maya diagram $\delta=(1^{\infty},0,1,0,1,0^{\infty})$ associated with $T$ and $r=4$. In this case, we have $W(\delta)=\{1,3,5\}$ and $R(\delta)=\{2,4\}$.}\label{Fig:Delta}
\end{figure}

If we write $W(\delta)=\{a_1<a_2<\cdots<a_{k+1}\}$ and $R(\delta)=\{b_1<b_2<\cdots<b_k\}$ for some $k\in\bb{Z}_{\geq0}$, then we have
\begin{align*}
	a_1<b_1<a_2<b_2<\cdots<a_k<b_k<a_{k+1}.
\end{align*}
In particular, the two sets $\{a\in W(\delta)\mid a>c\}$ and $\{b\in R(\delta)\mid b>c\}$ have the same cardinality, and hence $\varphi_c(\delta;q)>0$ for any $c\in W(\delta)$ and $q\in\bb{R}_{>0}$. Moreover, we have
\begin{align*}
	\sum_{c\in W(\delta)}\varphi_c(\delta;q)=1
\end{align*}
by \cite[Theorem~2]{Ker93}\footnote{When $q=1$, this kind of identity seems to be first explicitly stated and proved by Vershik \cite{Ver89}. Kirillov \cite{Kir89} gave another proof of this identity at $q=1$ and stated this $q$-analogue. See also Remark~\ref{Rem:Res}.}, and hence one may use these quantities to define transition probabilities of a Markov chain.

\begin{remark}
One may consider $\varphi_c(\delta;q)$ as a transition probability from $\delta$ to $\delta'$ obtained by exchanging the $(c-1)$-th and $c$-th components of $\delta$. By the well-known bijection between $\msf{Maya}$ and the set of charged Young diagrams (see for example \cite{MJD00}), this process can be viewed as a growth process of (charged) Young diagrams. Interestingly, such a growth process at $q=1$ has been studied in a different context under the name \emph{Plancherel growth process} (see \cite{Ker03}). Its $q$-analogue has also been studied, for example, by \cite{DMN21,FM12,Str08}.
\end{remark}

\subsubsection{Probabilistic models}

Now we describe our formula for $\msf{P}^{(r)}_{n}$. First, we note that for any $T\in\msf{SYT}_{n-1}$ and $c\in W(\boldsymbol{\delta}^{(r)}(T))$, the height of the $(c-1)$-th column of $T$ must be strictly greater than the height of the $c$-th column of $T$. Indeed, otherwise the standardness of $T$ would imply $\msf{top}_T(c-1)<\msf{top}_T(c)\leq r$, contradicting the condition $\delta_{c-1}=1$. In particular, we may add $\boxed{n}$ on top of $c$-th column of $T$ to obtain a new standard Young tableau.

\begin{dfn}\label{Def:f_c}
Let $T\in\msf{SYT}_{n-1}$ and $c\in\bb{Z}_{\geq0}$. If one can add $\boxed{n}$ on top of the $c$-th column of $T$, then we denote the resulting standard Young tableau by $\mf{f}_c(T)\in\msf{SYT}_n$.
\end{dfn}

\begin{dfn}\label{Def:Stochastic_Process}
Let $\msf{e}:\bb{Z}_{>0}\rightarrow\bb{Z}$ be a map and $q\in\bb{R}_{>0}$. We define a time-inhomogeneous Markov chain $\{\widetilde{\msf{X}}^{(\msf{e})}_{n}\}_{n\in\bb{Z}_{\geq0}}$ with state space $\msf{SYT}$ as follows:
\begin{itemize}
	\item The initial state is the empty standard Young tableau, i.e., $\widetilde{\msf{X}}^{(\msf{e})}_{0}=\emptyset$.
	\item The transition probability from $T\in\msf{SYT}_{n-1}$ to $T'\in\msf{SYT}_{n}$ at step $n$ is given by $\msf{P}^{(\msf{e}(n))}_{n}(T;T')$, where we set
	\begin{align*}
		\msf{P}^{(r)}_{n}(T;T')\coloneqq\begin{cases}
			\varphi_{c}\bigl(\boldsymbol{\delta}^{(r)}(T);q\bigr) & \text{if }T'=\mf{f}_c(T) \text{ with }c\in W\bigl(\boldsymbol{\delta}^{(r)}(T)\bigr),\\
			0 & \text{otherwise}.
		\end{cases}
	\end{align*}
\end{itemize}
This defines a stochastic process $\{\msf{X}^{(\msf{e})}_n\}_{n\in\bb{Z}_{\geq0}}$ with state space $\msf{Par}$ by $\msf{X}^{(\msf{e})}_n\coloneqq\pi\circ\widetilde{\msf{X}}^{(\msf{e})}_{n}$, where $\pi:\msf{SYT}\rightarrow\msf{Par}$ is the forgetful map that sends a tableau to its shape.
\end{dfn}

We note that the probability distribution of $\{\msf{X}^{(\msf{e})}_n\}_{n\in\bb{Z}_{\geq0}}$ at time $n$ depends  only on the first $n$ components of $\msf{e}$.  In particular, a conjugate Hessenberg function $\msf{e}:[n]\rightarrow\bb{Z}$ defines a stochastic process $\{\msf{X}^{(\msf{e})}_i\}_{0\leq i\leq n}$. We also remark that while the probabilistic interpretation requires $q\in\bb{R}_{>0}$, the process can be viewed algebraically for general $q$ by treating the transition probabilities as rational functions. See Figure~\ref{Fig:Transition} for an example of the transition probabilities.

\begin{figure}
\begin{tikzpicture}
\node at (0,0) {
\begin{ytableau}[] 
4 & *(red) 5 \\
1  & 2 & 3  & *(red) 6
\end{ytableau}};
\node at (5,3) {
\begin{ytableau}[] 
*(green) 7 \\
4 & *(red) 5 \\
1  & 2 & 3  & *(red) 6
\end{ytableau}};
\node at (5,0) {
\begin{ytableau}[] 
4 & *(red) 5 & *(green) 7 \\
1  & 2 & 3  & *(red) 6
\end{ytableau}};
\node at (5.315,-3) {
\begin{ytableau}[] 
4 & *(red) 5  \\
1  & 2 & 3  & *(red) 6 & *(green) 7
\end{ytableau}};
\draw[->,thick] (1.5,0.5) --node[above]{$\mf{f}_1$} (3.5,3);
\draw[->,thick] (1.5,0) --node[above]{$\mf{f}_3$} (3.5,0);
\draw[->,thick] (1.5,-0.5) --node[above]{$\mf{f}_5$} (3.5,-3);
\node[font=\Large] at (9,3) {$\frac{[-1]_q[-3]_q}{[-2]_q[-4]_q}$};
\node[font=\Large] at (9,0) {$\frac{[1]_q[-1]_q}{[2]_q[-2]_q}$};
\node[font=\Large] at (9,-3) {$\frac{[3]_q[1]_q}{[4]_q[2]_q}$};
\end{tikzpicture}
\caption{An example of the transition probabilities at step 7 for $\msf{e}(7)=4$.}\label{Fig:Transition}
\end{figure}

\subsubsection{Probabilistic formula}

We are now in a position to state our main theorem. 

\begin{thm}\label{Thm:Main}
	Let $\msf{e}:[n]\rightarrow\bb{Z}$ be a conjugate Hessenberg function and $\Gamma$ the corresponding unit interval graph. Then we have
\begin{align*}
	\mb{X}_{\Gamma}(x;q)=\sum_{\lambda\vdash n}\Prob\left\{\msf{X}^{(\msf{e})}_n=\lambda\right\}\cdot\prod_{i}[\lambda_i]_q!\cdot e_{\lambda}(x),
\end{align*}
where the probabilities are viewed as rational functions in $q$.
\end{thm} 

\begin{corollary}\label{Cor:qSS}
For any unit interval graph $\Gamma$, the $e$-expansion coefficients $c_{\lambda}(\Gamma;q)$ of the chromatic quasisymmetric function $\mb{X}_{\Gamma}(x;q)$ evaluated at $q\in\bb{R}_{>0}$ are nonnegative. 
\end{corollary}

By the result of Guay-Paquet \cite{GP13} (see Theorem~\ref{Thm:GP}), the Stanley--Stembridge conjecture has been reduced to the case of unit interval graphs. Hence Corollary~\ref{Cor:qSS} in particular proves the Stanley--Stembridge conjecture (see Conjecture~\ref{Conj_SS}) affirmatively.

\begin{corollary}
The Stanley--Stembridge conjecture is true.
\end{corollary}

\subsection{Organization of the paper}

This paper is organized as follows.

Before proceeding to the proof of the main theorem, we dedicate Section~\ref{Sec:Survey} to a detailed survey of the Stanley--Stembridge conjecture and related topics. In addition to fixing the notation used in our proof, there are two primary motivations for this exposition. First, we aim to illustrate the richness of the conjecture, which contrasts with the simplicity of its statement. In particular, we highlight how graph theory, combinatorics, geometry, and representation theory are interconnected, with Catalan objects serving as a common thread. Second, by demonstrating that the chromatic quasisymmetric functions of unit interval graphs---the main focus of this paper---admit generalizations in various directions, we intend to suggest potential future research directions. Readers interested primarily in the proof of the main theorem may safely skip this section after checking the essential notation.

Our proof of the main theorem relies on the characterization of the chromatic quasisymmetric functions of unit interval graphs in terms of the modular law (see Section~\ref{Sec:Modular_Law} for more details) given by Abreu and Nigro \cite{AN21}. The modular law relates values of a function on the set of conjugate Hessenberg functions, connecting local modifications of a conjugate Hessenberg function subject to certain non-local conditions. In Section~\ref{Sec:Reduction}, we first restate our main result in linear algebraic terms and then reduce it to three key lemmas, which essentially concern the local parts of the modular law. This section also aims to clarify how the non-local conditions are applied within the proof. 

Section~\ref{Sec:Key} is devoted to the proof of the key lemmas. Our basic tool is Corollary~\ref{Cor:Varphi_Useful} which relates $\varphi_c(\delta;q)$ and $\varphi_c(\delta';q)$, where $\delta'$ is obtained by modifying one component of $\delta$. The proofs are conceptually straightforward applications of Corollary~\ref{Cor:Varphi_Useful}, but in order to avoid a lengthy case-by-case analysis, we also introduce a rational function associated with a Maya diagram $\delta$, whose residues coincide with the values $\varphi_c(\delta;q)$. 

Finally, in Appendix~\ref{Sec:Prob}, we recall the basic language of probability theory for the convenience of the reader.

\subsection{Recent developments}

Since the first version of this paper appeared on the arXiv, there have been several further developments in this area. 

One of the obvious remaining problems is to prove the refined Stanley--Stembridge conjecture (Conjecture~\ref{Conj_SW}). Although our main result provides a formula for the $e$-expansion of the chromatic quasisymmetric functions of unit interval graphs, it is a sum of genuine rational functions in $q$. Thus, this result does not directly imply the refined conjecture. To approach the refined conjecture via our probabilistic formula, one needs to understand how the denominators cancel out after summing over all contributions.

Partial progress in this direction was made by Tom and Vailaya \cite{TV25}, who showed that summing over certain partial contributions yields polynomials in $q$. The resulting symmetric functions were identified by Huh, Hwang, Kim, Kim, and Oh in \cite{HHKKO25} with the $g$ functions defined by Abreu and Nigro \cite{AN23}. Consequently, they proved the $e$-positivity of these functions at $q=1$, refining our main result. 

Another refinement of our result was obtained by Griffin, Mellit, Romero, Weigl, and Wen \cite{GMRWW25}. They provided a formula for the expansion of chromatic quasisymmetric functions of unit interval graphs in terms of plethystically evaluated modified Macdonald polynomials, whose specialization at $t=1$ yields the $e$-expansion. Their proof utilizes the theory of $\bb{A}_{q,t}$ algebra, developed by Carlsson and Mellit \cite{CM18} in their proof of the shuffle conjecture, and also provides an alternative proof of the Stanley--Stembridge conjecture.

\subsection*{Acknowledgements}

The author thanks Syu Kato for introducing him to the subject and for helpful discussions. The author is also grateful to Mathieu Guay-Paquet for suggesting a reformulation of the transition probabilities, which clarified several proofs. 

The author also thanks the organizers and the other speakers of the conference \emph{Combinatorics on flag varieties and related topics 2025} held at Ajou University, where the author learned most of the material presented in Section~\ref{Sec:Survey}. Finally, the author thanks the anonymous referee for numerous constructive suggestions.

\section{Survey}\label{Sec:Survey}

In this section, we provide a survey of the Stanley--Stembridge conjecture and related topics. We also review various combinatorial objects related to unit interval graphs. As this survey introduces a variety of notation specific to this section, the notation and results to be used throughout the paper are emphasized by presenting them in theorem-like environments (e.g., Definition~\ref{Dfn_conjHess}) or are summarized in Section~2.1. 

\subsection{Notation}

For $n\in\bb{Z}_{>0}$, we write $[n]\coloneqq\{1,2,\ldots,n\}$. For $m\in\bb{Z}$, the $q$-integer $[m]_q$ is defined by 
\begin{align*}
	[m]_q\coloneqq\frac{1-q^m}{1-q},
\end{align*}
and for $n\in\bb{Z}_{>0}$, the $q$-factorial is defined by $[n]_q!\coloneqq[n]_q[n-1]_q\cdots[1]_q$.

\subsubsection{Symmetric functions}

We assume that the readers are familiar with the theory of symmetric functions and follow the basic notation in \cite{Mac95}. We use the power sum symmetric functions $p_{\lambda}$, the elementary symmetric functions $e_{\lambda}$, the complete symmetric functions $h_{\lambda}$, the monomial symmetric functions $m_{\lambda}$, and the Schur functions $s_{\lambda}$. We denote by $\Lambda$ the ring of symmetric functions over $\bb{Z}$ and by $\Lambda_n\subset\Lambda$ the submodule consisting of symmetric functions of degree $n$. There is an algebra involution $\omega$ on $\Lambda$ such that $\omega(h_\lambda)=e_{\lambda}$ and $\omega(s_{\lambda})=s_{\lambda^t}$, where $\lambda^t$ is the transpose of $\lambda$.

Let $\{b_{\lambda}\}_{\lambda}$ be a basis of the ring of symmetric functions indexed by partitions. We say that a symmetric function is \emph{$b$-positive} if the coefficients of its expansion in the basis $\{b_{\lambda}\}_{\lambda}$ are all nonnegative.

\subsubsection{Young tableaux}

In this paper, we adopt the French convention for Young diagrams, as it naturally aligns with the intuition of a stochastic process where boxes are dropped one by one. The Young diagram of a partition $\lambda=(\lambda_1,\lambda_2,\ldots)$ is defined as the set of boxes $\{(i,j)\in\bb{Z}_{>0}^2\mid 1\leq i\leq \lambda_{j}\}$. 
For a skew Young diagram $\nu$, we denote by $\msf{SSYT}(\nu)$ (resp. $\msf{SYT}(\nu)$) the set of \emph{semistandard} (resp. \emph{standard}) \emph{Young tableaux} of shape $\nu$. These are maps $T:\nu\rightarrow\bb{Z}_{>0}$ which are weakly increasing in each row and strictly increasing in each column (resp. and are also bijections from $\nu$ onto $\{1,2,\ldots,|\nu|\}$). 

\subsubsection{Graphs and posets}

A \emph{(simple) graph} $\Gamma=(V,E)$ consists of a set of vertices $V$ and a set of edges $E$, where each edge is a 2-element subset of $V$. By abuse of notation, we write $i\in\Gamma$ for a vertex $i\in V$ and $\{i,j\}\in\Gamma$ for an edge $\{i,j\}\in E$. 

A poset is called $(a+b)$-\emph{free} if it contains no induced subposet isomorphic to the disjoint union of an $a$-element chain and a $b$-element chain. The \emph{incomparability graph} of a poset $P$ is the graph whose vertex set is the set of elements of $P$, and where two distinct elements are connected by an edge if and only if they are incomparable in $P$.

\subsection{Stanley--Stembridge conjectures}

We first review the original formulation of the Stanley--Stembridge conjecture \cite{SS93} which originated in the study of the immanants of the Jacobi--Trudi matrices. 

\subsubsection{Characters of symmetric groups}

It is well-known (see, e.g., \cite{Mac95}) that the irreducible representations of the symmetric group $\mf{S}_n$ of degree $n$ are indexed by partitions of $n$. Let $\chi^{\lambda}$ denote the character of the irreducible representation of $\mf{S}_n$ corresponding to a partition $\lambda\vdash n$. For example, $\chi^{(n)}$ is the trivial character and $\chi^{(1^n)}$ is the sign character. 

Let $\msf{R}_n$ be the $\bb{Z}$-module generated by $\{\chi^{\lambda}\}_{\lambda\vdash n}$, i.e., the $\bb{Z}$-module of characters of virtual representations of $\mf{S}_n$. The \emph{Frobenius characteristic} $\msf{ch}:\msf{R}_n\rightarrow\Lambda_n$ is the $\bb{Z}$-module isomorphism defined by
\begin{align*}
	\msf{ch}(\chi)\coloneqq\frac{1}{n!}\sum_{w\in\mf{S}_n}\chi(w)p_{\rho(w)},
\end{align*}
where $\rho(w)\vdash n$ is the cycle type of $w\in\mf{S}_n$. For example, we have $\msf{ch}(\chi^{\lambda})=s_{\lambda}$. We define the \emph{monomial character} $\phi^{\lambda}\in\msf{R}_n$ by the condition $\msf{ch}(\phi^{\lambda})=m_{\lambda}$.

\subsubsection{Immanants of the Jacobi--Trudi matrices}

The immanant of a matrix was introduced by Littlewood and Richardson \cite{LR34}. For an $n\times n$ matrix $A=(a_{ij})$ and $\chi\in\msf{R}_n$, the \emph{immanant} of $A$ with respect to $\chi$ is defined by
\begin{align*}
	\msf{Imm}_{\chi}(A)\coloneqq\sum_{w\in\mf{S}_n}\chi(w)\prod_{i=1}^{n}a_{i,w(i)}.
\end{align*}
For example, $\msf{Imm}_{\chi^{(1^n)}}(A)$ is the determinant $\det A$, and $\msf{Imm}_{\chi^{(n)}}(A)$ is the \emph{permanent} of $A$. 

The \emph{Jacobi--Trudi matrix} $H_{\mu/\nu}(x)$ associated with a skew Young diagram $\mu/\nu$ satisfying $\ell(\mu)\leq n$ is defined by
\begin{align*}
	H_{\mu/\nu}(x)=\left(h_{\mu_i-\nu_j-i+j}(x)\right)_{1\leq i,j\leq n},
\end{align*}
where we understand that $h_0=1$ and $h_{r}=0$ if $r<0$. The \emph{Jacobi--Trudi formula} states that $\det H_{\mu/\nu}(x)$ is given by the skew Schur function $s_{\mu/\nu}(x)$. 

In \cite{GJ92}, Goulden and Jackson initiated a study of immanants of Jacobi--Trudi matrices and conjectured that $\msf{Imm}_{\chi^{\lambda}}H_{\mu/\nu}(x)$ is $m$-positive. This conjecture was proved by Greene \cite{Gre92}. In \cite{Ste92}, Stembridge proposed several refinements of this conjecture. One of them states that $\msf{Imm}_{\chi^{\lambda}}H_{\mu/\nu}(x)$ is $s$-positive, which was proved by Haiman \cite{Hai93}. The strongest conjecture in \cite{Ste92} asserts that $\msf{Imm}_{\phi^{\lambda}}H_{\mu/\nu}(x)$ is $s$-positive. 

To attack this conjecture, Stanley and Stembridge \cite{SS93} introduced a symmetric function $E^{\theta}_{\mu/\nu}(x)$ for $\theta\vdash N\coloneqq|\mu|-|\nu|$, defined by 
\begin{align*}
	\sum_{\lambda\vdash n}s_{\lambda}(x)\msf{Imm}_{\chi^{\lambda}}H_{\mu/\nu}(y)=\sum_{\theta\vdash N}E^{\theta}_{\mu/\nu}
(x)s_{\theta}(y).
\end{align*}
Since we have
\begin{align*}
\sum_{\lambda\vdash n}s_{\lambda}(x)\msf{Imm}_{\chi^{\lambda}}H_{\mu/\nu}(y)=\sum_{\lambda\vdash n}h_{\lambda}(x)\msf{Imm}_{\phi^{\lambda}}H_{\mu/\nu}(y),
\end{align*}
the $s$-positivity of $\msf{Imm}_{\phi^{\lambda}}H_{\mu/\nu}(x)$ for all $\lambda\vdash n$ is equivalent to the $h$-positivity of $E^{\theta}_{\mu/\nu}(x)$ for all $\theta\vdash N$. They conjectured \cite[Conjecture 5.1]{SS93} that $E^{\theta}_{\mu/\nu}(x)$ can be written as a nonnegative linear combination of functions of the form $E^{(N)}_{\mu'/\nu'}(x)$. This conjecture remains open (see \cite{Les24} for a partial result), but it would reduce the $h$-positivity of $E^{\theta}_{\mu/\nu}(x)$ to that of $E^{(N)}_{\mu/\nu}(x)$.

Concretely, $E^{\theta}_{\mu/\nu}(x)$ is given by
\begin{align*}
	E^{\theta}_{\mu/\nu}(x)=\sum_{w\in\mf{S}_n}K_{\theta,\mu+\delta-w(\nu+\delta)}p_{\rho(w)}(x),
\end{align*}
where $K_{\theta,\alpha}$ is the \emph{Kostka number}, i.e., the number of semistandard Young tableau $T$ of shape $\theta$ and weight $\alpha$. Here, we understand that $K_{\theta,\alpha}=0$ if $\alpha_i<0$ for some $i$. Since $K_{(N),\alpha}=1$ for any composition $\alpha\models N$, $E^{(N)}_{\mu/\nu}(x)$ depends only on  $\sigma\coloneqq\{(i,j)\in[n]^2\mid \mu_i-\nu_j-i+j<0\}$ which induces a partial order $<_{\sigma}$ on $[n]$ by
\begin{align*}
	j<_{\sigma}i\iff (i,j)\in\sigma.
\end{align*}

In fact, Stanley--Stembridge \cite{SS93} showed that $E^{(N)}_{\mu/\nu}(x)$ can be identified with $\omega\mb{X}_{P_{\sigma}}(x)$ in terms of the poset $P_{\sigma}\coloneqq([n],<_{\sigma})$, where we set
\begin{align*}
	\mb{X}_{P}(x)\coloneqq\sum_{\lambda\vdash n}\#\Bigl\{\text{partitions of }P\text{ into ordered chains of cardinality }\lambda_1,\lambda_2,\ldots\Bigr\}m_{\lambda}(x)
\end{align*}
for a poset $P$ of size $n$. Hence it is a natural problem to consider when $\mb{X}_{P}(x)$ is $e$-positive. 

Since posets of the form $P_{\sigma}$ can be characterized as both $(\mb{3}+\mb{1})$-free and $(\mb{2}+\mb{2})$-free posets by the result of \cite{DK68}, they studied how the $(\mb{3}+\mb{1})$-free and $(\mb{2}+\mb{2})$-free conditions affect the $e$-positivity of $\mb{X}_{P}(x)$. Based on numerical experiments, they conjectured the following, which is commonly referred to as \emph{the Stanley--Stembridge conjecture}.

\begin{conj}[Stanley--Stembridge {\cite[Conjecture 5.5]{SS93}}]\label{Conj_SS}
For any $(\mb{3}+\mb{1})$-free poset $P$, $\mb{X}_{P}(x)$ is $e$-positive.
\end{conj}

\subsubsection{Chromatic symmetric functions}

In his study of the four color problem, Birkhoff \cite{Bir12} introduced the notion of the chromatic polynomial for graphs originally in the context of planar graphs. This polynomial counts the number of \emph{proper colorings}, i.e., colorings of the vertices such that two vertices have different colors if they are connected by an edge. 

In 1995, Stanley \cite{Sta95} introduced the notion of the chromatic symmetric function $\mb{X}_{\Gamma}(x)$ for a graph $\Gamma$, which contains more refined information on the number of proper colorings with $\lambda_i$ vertices colored by $i\in\bb{Z}_{>0}$ for any partition $\lambda=(\lambda_1,\lambda_2,\ldots)$, i.e., 
\begin{align*}
	\mb{X}_{\Gamma}(x)\coloneqq\sum_{\lambda\vdash n}\#\Bigl\{\text{proper colorings of }\Gamma\text{ with }\lambda_i\text{ vertices colored by }i\Bigr\}m_{\lambda}(x),
\end{align*}
where $n$ is the number of vertices of $\Gamma$. Stanley \cite{Sta95} studied basic properties of chromatic symmetric functions, such as their expansions in terms of power sum and elementary symmetric functions.

Since $\mb{X}_{\Gamma}(x)$ coincides with $\mb{X}_{P}(x)$ when $\Gamma$ is the incomparability graph of a poset $P$, the Stanley--Stembridge conjecture is equivalent to the following statement.

\begingroup 
    \renewcommand{\thethm}{\ref{Conj_SS}'} 
    \addtocounter{thm}{-1} 
    \begin{conj}[{\cite[Conjecture 5.1]{Sta95}}]
    For any incomparability graph $\Gamma$ of a $(\bf 3+1)$-free poset, $\mb{X}_{\Gamma}(x)$ is $e$-positive.
    \end{conj}
\endgroup

There have been many studies on the $e$-positivity of chromatic symmetric functions for graphs not necessarily coming from $(\bf 3+1)$-free posets. For example,  \cite{AWW24,BCCCGKKLLS25,CHW26,CH22,CH19,Dah18,DFW20,DSW20,DW18,FHK19,FKKMMT21,GS01,HHT19,HNY20, LLWY21,LY21,MPW24,MS25,QTW24,Tom25,Tsu18,Wan24,WW23a,WW23b,WZ25,Wan22,Zhe22} prove the $e$-positivity (or non-$e$-positivity) for special classes of graphs. A weaker statement regarding the $s$-positivity for any $(\bf 3+1)$-free poset was proved by Gasharov \cite{Gas96}, by exhibiting a combinatorial interpretation of the coefficients in the Schur function expansion. 

\subsection{Unit interval graphs}\label{Sec:UIG}

A significant breakthrough relevant to this paper was made by Guay-Paquet \cite{GP13}, who proved that it is sufficient to verify the $e$-positivity for posets that are both $(\bf 3+1)$-free and $(\bf 2+2)$-free to establish the Stanley--Stembridge conjecture. 

\begin{thm}[\cite{GP13}]\label{Thm:GP}
If $\mb{X}_{P}(x)$ is $e$-positive for all posets $P$ which are both $(\bf 3+1)$-free and $(\bf 2+2)$-free, then $\mb{X}_{P}(x)$ is $e$-positive for all $(\bf 3+1)$-free posets $P$.
\end{thm}

It is well known \cite{SS58} that a poset is $(\bf 3+1)$-free and $(\bf 2+2)$-free if and only if it is isomorphic to a \emph{unit interval order}, i.e., a poset whose underlying set is a collection of unit intervals on the real line, partially ordered by the relation $[a,a+1]<[b,b+1]$ if and only if $a+1<b$. It is natural to label the intervals with $[n]\coloneqq\{1,2,\ldots,n\}$ in order from left to right. With this labeling, a poset $([n],<_P)$ is isomorphic to a unit interval order if and only if it satisfies the following conditions \cite[Proposition 4.1]{SW16}:
\begin{itemize}
	\item $i<_Pj$ implies $i<j$ in the natural order, 
	\item if $i,j,k\in[n]$ satisfy $i<_Pk$ and $j$ is incomparable with both $i$ and $k$, then we have $i<j<k$ in the natural order.
\end{itemize}

The incomparability graph of such a poset is known as a \emph{(natural) unit interval graph}. Since the number of unit interval graphs with $n$ vertices coincides with the $n$-th Catalan number $C_n=\frac{1}{n+1}\binom{2n}{n}$, which enumerates various combinatorial objects (see for example \cite{Sta15}), Guay-Paquet's result allows us to approach the Stanley--Stembridge conjecture using various combinatorial objects, such as Dyck paths, Hessenberg functions, area sequences, and 312-avoiding permutations. We summarize several bijections between these objects below.

\subsubsection{Dyck paths}

A lattice path from $(0,0)$ to $(n,n)$ consisting of north steps $N=(0,1)$ and east steps $E=(1,0)$ is called a \emph{Dyck path} of length $n$ if it never goes under the diagonal $y=x$. It is well known that the number of Dyck paths of length $n$ is given by the Catalan number $C_n$. 

There are several natural ways to construct a unit interval graph from a Dyck path $\pi$. One approach is to use the diagonal boxes as the set of vertices. We label these boxes with $[n]$ from bottom to top. Two diagonal boxes are connected by an edge if they can be connected by a path which first travels to the north and then to the east while remaining below $\pi$. In other words, $i<j\in[n]$ are connected by an edge if and only if the box with top right corner $(i,j)$ lies under the path $\pi$. Let $\Gamma_{\pi}$ denote the resulting graph. The map $\pi\mapsto\Gamma_{\pi}$ defines a bijection between the set of Dyck paths and the set of unit interval graphs. 

Another way is to use the boxes immediately to the right of the north steps as the set of vertices. We label these vertices with $[n]$ in the \emph{reading order}, i.e., we read the boxes diagonal by diagonal, from bottom to top, and from left to right within each diagonal. We connect two vertices $u<v\in[n]$ if and only if either
\begin{itemize}
	\item $c(v)=c(u)$ and $v$ is above $u$, or
	\item $c(v)=c(u)+1$ and $v$ is below $u$.
\end{itemize}
Here, we define the \emph{content} $c(u)$ of a box $u=(i,j)$ by $c(u)=j-i$. Let $\Gamma'_{\pi}$ denote the resulting graph. This correspondence $\pi\mapsto\Gamma'_{\pi}$ also yields a bijection between the set of Dyck paths and the set of unit interval graphs. In fact, one can construct a bijection $\zeta$ on the set of Dyck paths such that $\Gamma'_{\pi}=\Gamma_{\zeta(\pi)}$, by the proof of \cite[Theorem 3.15]{Hag08} (see also \cite{CM18}). This correspondence is useful for describing the connection with unicellular LLT polynomials (see Section~\ref{Sec:LLT}). Figure~\ref{Fig_Dyck_UIG} illustrates the two constructions of unit interval graphs.

\begin{figure}
\centering
\begin{tikzpicture}
\node (Dyck) at (0,0){
\begin{tikzpicture}[scale=0.5,edge_style/.style={thick, rounded corners=8pt}]
\foreach \x in {0,...,5} \draw[gray!60] (\x,0) -- (\x,5);
\foreach \y in {0,...,5} \draw[gray!60] (0,\y) -- (5,\y);
\draw[gray!60] (0,0) -- (5,5);
\draw[->,very thick] (0,0)--(0,1)--(0,2)--(0,3)--(1,3)--(1,4)--(2,4)--(3,4)--(3,5)--(4,5)--(5,5);
\tikzset{cellnode/.style={minimum width=1, minimum height=1}} 
    \node[cellnode] (N1) at (0.5, 0.5) {1};
    \node[cellnode] (N2) at (1.5, 1.5) {2};
    \node[cellnode] (N3) at (2.5, 2.5) {3};
    \node[cellnode] (N4) at (3.5, 3.5) {4};
    \node[cellnode] (N5) at (4.5, 4.5) {5};
    \draw[edge_style,red] (N1.north) -- (N1.north |- N2.west) -- (N2.west);
    \draw[edge_style,red] (N1.north) -- (N1.north |- N3.west) -- (N3.west);
    \draw[edge_style,red] (N2.north) -- (N2.north |- N3.west) -- (N3.west);
    \draw[edge_style,red] (N2.north) -- (N2.north |- N4.west) -- (N4.west);
    \draw[edge_style,red] (N3.north) -- (N3.north |- N4.west) -- (N4.west) ;
    \draw[edge_style,red] (N4.north) -- (N4.north |- N5.west) -- (N5.west);
    \node at (2.5,-1) {$\pi=NNNENEENEE$};
\end{tikzpicture}};
\node[anchor=west] (hgraph) at (2.5,0) {
\begin{tikzpicture}[scale=0.5]
\fill (0,0) circle (5pt);
\node[below=3pt] at (0,0) {1};
\fill (2,0) circle (5pt);
\node[below=3pt] at (2,0) {2};
\fill (4,0) circle (5pt);
\node[below=3pt] at (4,0) {3};
\fill (6,0) circle (5pt);
\node[below=3pt] at (6,0) {4};
\fill (8,0) circle (5pt);
\node[below=3pt] at (8,0) {5};
\draw (0,0)--(2,0);
\draw (2,0)--(4,0);
\draw (4,0)--(6,0);
\draw (6,0)--(8,0);
\draw (0,0) to[out=60,in=120] (4,0);
\draw (2,0) to[out=60,in=120] (6,0);
\node at (2,-2.5) {$\Gamma_{\pi}=\Gamma'_{\pi'}$};
\end{tikzpicture}};
\node (Dyck) at (10,0){
\begin{tikzpicture}[scale=0.5,node_style/.style={minimum width=1, minimum height=1}]\foreach \x in {0,...,5} \draw[gray!60] (\x,0) -- (\x,5);
\foreach \y in {0,...,5} \draw[gray!60] (0,\y) -- (5,\y);
\draw[gray!60] (0,0) -- (5,5);
\draw[->,very thick] (0,0)--(0,1)--(0,2)--(1,2)--(2,2)--(2,3)--(3,3)--(3,4)--(3,5)--(4,5)--(5,5);
    \node[node_style] (N1) at (0.5, 0.5) {1};
    \node[node_style] (N4) at (0.5, 1.5) {4};
    \node[node_style] (N2) at (2.5, 2.5) {2};
    \node[node_style] (N3) at (3.5, 3.5) {3};
    \node[node_style] (N5) at (3.5, 4.5) {5};
    \draw[thick,red] (N1) to[bend right=30] (N2);
    \draw[thick,red] (N2) to[bend right=30] (N3);
    \draw[thick,red] (N1) to[bend right=50] (N3);
    \draw[thick,red] (N2) to[bend right=30] (N4);
    \draw[thick,red] (N3) to[bend right=40] (N4);
    \draw[thick,red] (N4) to[bend left=50] (N5);
    \node at (2.5,-1) {$\pi'=NNEENENNEE$};
\end{tikzpicture}};
\end{tikzpicture}
\caption{Dyck paths and corresponding unit interval graphs}\label{Fig_Dyck_UIG}
\end{figure}

\subsubsection{Hessenberg functions}

A function $\msf{h}:[n]\rightarrow\bb{Z}$ is called a \emph{Hessenberg function} if it satisfies
\begin{align*}
	\msf{h}(i)\geq i,\hspace{1em}\msf{h}(i)\leq\msf{h}(i+1)
\end{align*}
for any $i\in[n]$, with the convention that $\msf{h}(n+1)=n$. For a Hessenberg function $\msf{h}$, we define a partial order $<_{\msf{h}}$ on $[n]$ given by
\begin{align*}
	i<_{\msf{h}}j\iff \msf{h}(i)<j.
\end{align*}
This correspondence induces a bijection between the set of Hessenberg functions and the set of isomorphism classes of unit interval orders (see \cite[Proposition 4.1]{SW16}). The incomparability graph of the poset $([n],<_{\msf{h}})$ is the unit interval graph $\Gamma_{\msf{h}}$ with vertex set $[n]$ and edge set $\{\{i,j\}\mid i<j\leq\msf{h}(i)\}$.

The bijection $\pi\mapsto\msf{h}$ between the set of Dyck paths and the set of Hessenberg functions, satisfying $\Gamma_{\pi}=\Gamma_{\msf{h}}$, is given by 
\begin{align*}
	\msf{h}(i)=\text{number of boxes in the }i\text{-th column under }\pi.
\end{align*}

\subsubsection{Conjugate Hessenberg functions}\label{Sec:conj_Hess}

In this paper, we use a variant of the Hessenberg function, which we call a conjugate Hessenberg function, defined as follows. 

\begin{dfn}\label{Dfn_conjHess}
	A function $\msf{e}:[n]\rightarrow\bb{Z}$ is called a \textit{conjugate Hessenberg function} if it satisfies
\begin{align*}
	0\leq\msf{e}(i)<i,\hspace{1em}\msf{e}(i)\leq\msf{e}(i+1)
\end{align*}
for any $i\in[n]$, with the convention that $\msf{e}(n+1)=\infty$. We denote by $\bb{E}_n$ the set of conjugate Hessenberg functions on $[n]$. For $\msf{e}\in\bb{E}_n$, we associate the unit interval graph $\Gamma_{\msf{e}}$ with vertex set $[n]$ and edge set $\left\{\{i,j\}\mid \msf{e}(j)<i<j\right\}$.
\end{dfn}

The bijection $\pi\mapsto\msf{e}$ between the set of Dyck paths and the set of conjugate Hessenberg functions satisfying $\Gamma_{\pi}=\Gamma_{\msf{e}}$ is given by 
\begin{align*}
	\msf{e}(i)=\text{number of boxes in the }i\text{-th row to the left of }\pi.
\end{align*}
One of the reasons for using conjugate Hessenberg functions is that they are well-suited for describing the inductive structure behind chromatic symmetric functions. 

\begin{figure}
\centering
\begin{tikzpicture}
\node[anchor=west] (graph) at (-4,0) {
\begin{tikzpicture}[scale=0.5]
\fill (0,0) circle (5pt);
\node[below=3pt] at (0,0) {1};
\fill (2,0) circle (5pt);
\node[below=3pt] at (2,0) {2};
\fill (4,0) circle (5pt);
\node[below=3pt] at (4,0) {3};
\fill (6,0) circle (5pt);
\node[below=3pt] at (6,0) {4};
\fill (8,0) circle (5pt);
\node[below=3pt] at (8,0) {5};
\draw (0,0)--(2,0);
\draw (2,0)--(4,0);
\draw (4,0)--(6,0);
\draw (6,0)--(8,0);
\draw (0,0) to[out=60,in=120] (4,0);
\draw (2,0) to[out=60,in=120] (6,0);
\node at (3.5,-2) {$\Gamma$};
\end{tikzpicture}};
\node (Dyck) at (4,0){
\begin{tikzpicture}[scale=0.5,edge_style/.style={thick, rounded corners=8pt}]
\foreach \x in {0,...,5} \draw[gray!60] (\x,0) -- (\x,5);
\foreach \y in {0,...,5} \draw[gray!60] (0,\y) -- (5,\y);
\draw[gray!60] (0,0) -- (5,5);
\draw[->,very thick] (0,0)--(0,1)--(0,2)--(0,3)--(1,3)--(1,4)--(2,4)--(3,4)--(3,5)--(4,5)--(5,5);
\node at (2.5,-1) {$\pi$};
\draw[<->,red] (0.5,0) -- (0.5,3);
\draw[<->,red] (1.5,0) -- (1.5,4);
\draw[<->,red] (2.5,0) -- (2.5,4);
\draw[<->,red] (3.5,0) -- (3.5,5);
\draw[<->,red] (4.5,0) -- (4.5,5);
\draw[<->,green] (0,3.5) -- (1,3.5);
\draw[<->,green] (0,4.5) -- (3,4.5);
\end{tikzpicture}};
\draw[<->,thick] ($(Dyck.west)+(0,0.3)$)--($(graph.east)+(0,0.3)$);
\node (h) at (8.5,1.5) {$\msf{h}=(3,4,4,5,5)$};
\node (e) at (8.5,-1) {$\msf{e}=(0,0,0,1,3)$};
\draw[<->,thick,red] (Dyck.east)+(0,0.4) -- (h.west);
\draw[<->,thick,green] (Dyck.east)+(0,0.1) -- (e.west);
\end{tikzpicture}
\caption{Bijections between Catalan objects: Dyck path $\pi$, Hessenberg function $\msf{h}=(3,4,4,5,5)$, conjugate Hessenberg function $\msf{e}=(0,0,0,1,3)$, and unit interval graph $\Gamma=\Gamma_{\pi}$. These objects also correspond to the area sequence $\msf{a}=(0,1,2,2,1)$ and the 312-avoiding permutation $w=34251$.}
\label{Fig_Catalan}
\end{figure}

\subsubsection{Area sequences}

We often view a unit interval graph as a directed graph by choosing a natural orientation for each edge. Namely, an edge between $i$ and $j$ is oriented from $i$ to $j$ whenever $i<j$. As a directed graph, a unit interval graph is a special case of a circular unit arc digraph, which is parametrized by an integer sequence called an \emph{area sequence} \cite{AP18}, i.e., a function $\msf{a}:[n]\rightarrow\bb{Z}$ satisfying 
\begin{align*}
	0\leq\msf{a}(i)<n,\hspace{1em}\msf{a}(i+1)\leq \msf{a}(i)+1
\end{align*}
for any $i\in[n]$, with the convention that $\msf{a}(n+1)=\msf{a}(1)$. The \emph{circular unit arc digraph} $\Gamma_{\msf{a}}$ associated with an area sequence $\msf{a}$ is the directed graph with vertex set $[n]$ and edges 
\begin{align*}
	(i-\msf{a}(i))\to i,\quad (i-\msf{a}(i)+1)\to i,\quad \ldots,\quad (i-1)\to i
\end{align*}
for all $i\in[n]$, where vertex indices are considered modulo $n$. The case $\msf{a}(1)=0$ corresponds to a unit interval graph. 

The bijection $\msf{a}\mapsto\msf{e}$ between the set of area sequences satisfying $\msf{a}(1)=0$ and the set of conjugate Hessenberg functions satisfying $\Gamma_{\msf{a}}=\Gamma_{\msf{e}}$ is given by
\begin{align*}
	\msf{e}(i)=i-1-\msf{a}(i)
\end{align*}
for all $i\in[n]$. In terms of the corresponding Dyck path $\pi$, $\msf{a}(i)$ is the number of boxes between $\pi$ and the diagonal in the $i$-th row from the bottom. 

\subsubsection{312-avoiding permutations}

A permutation $w=w_1w_2\cdots w_n$ (in one-line notation) is called \emph{312-avoiding} (or \emph{codominant}) if there exist no indices $i<j<k$ such that $w_i>w_k>w_j$. It is known \cite{Knu97} that the number of 312-avoiding permutations in $\mf{S}_n$ is also given by the Catalan number $C_n$. We consider the map from $\mf{S}_n$ to the set of Hessenberg functions defined by sending $w=w_1w_2\cdots w_n\in\mf{S}_n$ to the function $\msf{h}$ given by
\begin{align*}
\msf{h}(i)=\max\{w_1,w_2,\ldots,w_i\}
\end{align*}
for all $i\in[n]$. This map is surjective, and each fiber contains a unique 312-avoiding permutation (see \cite[Proposition~3.1]{Hai93}).

\subsection{Chromatic quasisymmetric functions}

Another important advance towards the Stanley--Stembridge conjecture was made by Shareshian and Wachs \cite{SW12,SW16}, who introduced a $q$-analogue $\mb{X}_{\Gamma}(x;q)$ of the chromatic symmetric function for any \emph{labeled graph} $\Gamma$, i.e., a graph with vertex set $[n]$ for some $n\in\bb{Z}_{>0}$. For a (not necessarily proper) coloring $\kappa:[n]\rightarrow\bb{Z}_{>0}$ of the vertices of $\Gamma$, we set
\begin{align}\label{Eqn_asc}
	\msf{asc}(\kappa)\coloneqq\#\left\{\{i,j\}\in\Gamma\mid i<j\mbox{ and }\kappa(i)<\kappa(j)\right\}.
\end{align}

\begin{dfn}[\cite{SW12,SW16}]\label{Dfn_CQF}
Let  $\Gamma$ be a labeled graph. The \emph{chromatic quasisymmetric function} $\mb{X}_{\Gamma}(x;q)$ of $\Gamma$ is defined by
\begin{align*}
	\mb{X}_{\Gamma}(x;q)\coloneqq\sum_{\substack{\kappa:\Gamma\rightarrow\bb{Z}_{>0}\\\kappa:\mbox{\tiny proper}}}q^{\msf{asc}(\kappa)}\prod_{i\in\Gamma}x_{\kappa(i)}.
\end{align*}
\end{dfn}

We summarize basic properties of chromatic quasisymmetric functions. 

\begin{prop}[\cite{SW16}]\label{Prop_basic_properties}
Let $\Gamma=([n],E)$ be a labeled graph.
\begin{enumerate}
    \item At $q=1$, $\mb{X}_{\Gamma}(x;q)$ reduces to the chromatic symmetric function $\mb{X}_{\Gamma}(x)$.

    \item If $\Gamma$ is the ordered disjoint union of two labeled graphs $\Gamma_1$ and $\Gamma_2$, then 
    \begin{align*}
    	\mb{X}_{\Gamma}(x;q) = \mb{X}_{\Gamma_1}(x;q) \mb{X}_{\Gamma_2}(x;q).
    \end{align*}
    \item The function $\mb{X}_{\Gamma}(x;q)$ is a quasisymmetric function in the sense of Gessel \cite{Ges84}.
    \item If $\Gamma$ is a unit interval graph with its natural labeling, then $\mb{X}_{\Gamma}(x;q)\in\Lambda[q]$.
    \item If $\Gamma$ is the complete graph $K_n$ on $n$ vertices, then
    \begin{align*}
        \mb{X}_{K_n}(x; q) = [n]_q! \, e_n(x).
    \end{align*}
    \item If $\mb{X}_{\Gamma}(x;q)\in\Lambda[q]$, then we have $\mb{X}_{\Gamma}(x;q)=q^{|E|}\mb{X}_{\Gamma}(x;q^{-1})$.
\end{enumerate}
\end{prop}

Recall that by the aforementioned work of Guay-Paquet \cite{GP13}, the Stanley--Stembridge conjecture was reduced to the $e$-positivity of $\mb{X}_{\Gamma}(x)$ for unit interval graphs $\Gamma$. Thus the following conjecture by Shareshian--Wachs refines the Stanley--Stembridge conjecture.

\begin{conj}[Shareshian--Wachs \cite{SW12,SW16}]\label{Conj_SW}
For any unit interval graph $\Gamma$, the $e$-expansion of $\mb{X}_{\Gamma}(x;q)$ has coefficients in $\bb{Z}_{\geq0}[q]$.
\end{conj}

In fact, the notion of chromatic quasisymmetric function was generalized by Ellzey \cite{Ell17} to any directed graph by replacing $i<j$ in (\ref{Eqn_asc}) with a directed edge $i\to j$. There are other classes of directed graphs whose chromatic quasisymmetric functions are symmetric. For example, the chromatic quasisymmetric functions for circular unit arc digraphs are symmetric by \cite{AP18,Ell17}. See also \cite{GPS25} on the problem of when chromatic quasisymmetric functions are symmetric. 

The $e$-positivity of chromatic quasisymmetric functions which are symmetric has also been studied by many authors, including \cite{AN21,AN23,AP18,CH19,CMP23,Ell17,HP19,HNY20,Hwa24,LS22,NT23,Tom25,Wan22}. For a unit interval graph $\Gamma$, $\mb{X}_{\Gamma}(x;q)$ has many incarnations in terms of other mathematical objects corresponding to $\Gamma$. In the following, we review its connections with LLT polynomials, Hessenberg varieties, and Hecke algebras.
   
\subsubsection{Unicellular LLT polynomials}\label{Sec:LLT}

The notion of the LLT polynomial was originally defined by Lascoux, Leclerc, and Thibon \cite{LLT97} in terms of ribbon tableaux. We briefly recall its definition as reformulated in \cite{HHL05,HHLRU05}. 

Let $\boldsymbol{\nu}=(\nu^{(1)},\ldots,\nu^{(k)})$ be a tuple of skew Young diagrams. Given a semistandard tableau $T=(T^{(1)},\ldots,T^{(k)})\in\msf{SSYT}(\boldsymbol{\nu})\coloneqq\msf{SSYT}(\nu^{(1)})\times\cdots\times\msf{SSYT}(\nu^{(k)})$, a pair of boxes $u\in\nu^{(i)}$ and $v\in\nu^{(j)}$ is called an \emph{attacking inversion} if $T^{(i)}(u)>T^{(j)}(v)$ and either
\begin{itemize}
	\item $i<j$ and $c(u)=c(v)$, or
	\item $i>j$ and $c(u)=c(v)+1$.
\end{itemize}
We denote by $\msf{inv}(T)$ the number of attacking inversions in $T$. The \emph{LLT polynomial} $\msf{LLT}_{\boldsymbol{\nu}}(x;q)$ associated with $\boldsymbol{\nu}$ is defined by
\begin{align*}
	\msf{LLT}_{\boldsymbol{\nu}}(x;q)\coloneqq\sum_{T\in\msf{SSYT}(\boldsymbol{\nu})}q^{\msf{inv}(T)}\prod_{i=1}^{k}\prod_{u\in\nu^{(i)}}x_{T^{(i)}(u)}.
\end{align*}
When each $\nu^{(i)}$ consists of a single box (resp. vertical strip), $\msf{LLT}_{\boldsymbol{\nu}}(x;q)$ is called a \emph{unicellular LLT polynomial} (resp. a \emph{vertical-strip LLT polynomial}). 

For an $n$-tuple of skew Young diagrams $\boldsymbol{\nu}=(\nu^{(1)},\ldots,\nu^{(n)})$ with $|\nu^{(i)}|=1$, we consider the unit interval order given by the collection of unit intervals 
\begin{align*}
	\left\{I_i\coloneqq\left[c(\nu^{(i)})-i\epsilon,c(\nu^{(i)})-i\epsilon+1\right]\mid i\in[n]\right\}
\end{align*}
for any sufficiently small $\epsilon>0$. Let $\Gamma_{\boldsymbol{\nu}}$ denote the corresponding unit interval graph. We note that any unit interval graph can be obtained as $\Gamma_{\boldsymbol{\nu}}$ for some $\boldsymbol{\nu}$. For example, one can construct $\boldsymbol{\nu}$ such that $\Gamma'_{\pi}=\Gamma_{\boldsymbol{\nu}}$ for a Dyck path $\pi$ by choosing $\nu^{(i)}$ to be the $i$-th box from the top in the set of vertices of $\Gamma'_{\pi}$.

The set $\msf{SSYT}(\boldsymbol{\nu})$ can be identified with the set of all colorings $\kappa:\Gamma_{\boldsymbol{\nu}}\rightarrow\bb{Z}_{>0}$ of the vertices of $\Gamma_{\boldsymbol{\nu}}$. A pair $(u,v)$ with $u\in\nu^{(i)}$, $v\in\nu^{(j)}$ forms an attacking inversion for $\kappa$ if and only if $I_j$ precedes $I_i$ in the natural order, $\kappa(I_j)<\kappa(I_i)$, and they are connected by an edge in $\Gamma_{\boldsymbol{\nu}}$. By comparing the definitions, one obtains 
\begin{align*}
	\msf{LLT}_{\boldsymbol{\nu}}(x;q)=\sum_{\kappa:\Gamma_{\boldsymbol{\nu}}\rightarrow\bb{Z}_{>0}}q^{\msf{asc}(\kappa)}\prod_{i\in\Gamma_{\boldsymbol{\nu}}}x_{\kappa(i)}.
\end{align*}
In other words, the unicellular LLT polynomial $\msf{LLT}_{\boldsymbol{\nu}}(x;q)$ can be defined in the same way as $\mb{X}_{\Gamma_{\boldsymbol{\nu}}}(x;q)$, except that the sum is over all  colorings. 

In fact, Carlsson and Mellit \cite[Proposition 3.5]{CM18} showed that 
\begin{align*}
	\mb{X}_{\Gamma_{\boldsymbol{\nu}}}(x;q)=(q-1)^{-n}\msf{LLT}_{\boldsymbol{\nu}}\left[(q-1)X;q\right].
\end{align*}
Here, the \emph{plethystic substitution} $f(x)\mapsto f[(q-1)X]$ for $f\in\Lambda\otimes\bb{Q}(q)$ is an algebra homomorphism $\Lambda\otimes\bb{Q}(q)\rightarrow\Lambda\otimes\bb{Q}(q)$ characterized by
\begin{align*}
	p_r(x)\mapsto (q^r-1)p_r(x).
\end{align*}

We note that as an analogue of the Stanley--Stembridge conjecture, the $e$-positivity of $\msf{LLT}_{\boldsymbol{\nu}}(x;q+1)$ for the vertical-strip LLT polynomial $\msf{LLT}_{\boldsymbol{\nu}}(x;q)$ has  also been studied by \cite{Ale21,AP18,AS22,DAd20,GHQR19}. 

\subsubsection{Regular semisimple Hessenberg varieties}\label{Sec:Hess}

The notion of a Hessenberg variety was introduced by De Mari, Procesi, and Shayman \cite{DPS92}. We first recall its general definition and then specialize to the type $A$ cases, which are related to the chromatic quasisymmetric functions of unit interval graphs.

Let $\msf{G}$ be a complex reductive algebraic group and $\msf{B}\subset \msf{G}$ be a Borel subgroup with Lie algebras $\mf{g}$ and $\mf{b}$. A subspace $\mf{m}\subset\mf{g}$ stable under the adjoint action of $\msf{B}$ defines a vector bundle $\msf{G}\times^{\msf{B}}\mf{m}\coloneqq(\msf{G}\times\mf{m})/\msf{B}$ over the flag variety $\msf{G}/\msf{B}$, where $\msf{B}$ acts on $\msf{G}\times\mf{m}$ by $b\cdot(g,m)=(gb^{-1},\msf{Ad}_{b}(m))$. In the case $\mf{m}=\mf{n}$ (the nilpotent radical of $\mf{b}$), the bundle $\msf{G}\times^{\msf{B}}\mf{n}$ can be identified with the cotangent bundle $T^{\ast}(\msf{G}/\msf{B})\cong\msf{G}\times^{\msf{B}}(\mf{g}/\mf{b})^{\ast}$ via an isomorphism $(\mf{g}/\mf{b})^{\ast}\cong\mf{n}$ induced by a $\msf{G}$-invariant nondegenerate bilinear form on $\mf{g}$. We define an analogue of the Springer resolution by 
\begin{align*}
	\begin{array}{rccc}
\mu_{\msf{m}}: & \msf{G}\times^{\msf{B}}\mf{m} & \longrightarrow & \mf{g} \\
   & \rotatebox[origin=c]{90}{$\in$} & & \rotatebox[origin=c]{90}{$\in$} \\
   & (g,m) & \longmapsto     & \msf{Ad}_{g}(m).
   \end{array}
\end{align*}
This morphism is projective due to the projectivity of $\msf{G}/\msf{B}$. When $\mf{m}=\mf{n}$, a non-empty fiber is called a \emph{Springer fiber}. When $\mf{m}$ contains $\mf{b}$, the morphism $\mu_{\mf{m}}$ is surjective, and its fiber is called a \emph{Hessenberg variety}. The fiber over a regular semisimple element is referred to as a \emph{regular semisimple Hessenberg variety}. It is known \cite{DPS92} that every regular semisimple Hessenberg variety is smooth of dimension $\dim\mf{m}/\mf{b}$. Consequently, the restriction $\mu^{\mr{rs}}_{\mf{m}}$ of $\mu_{\mf{m}}$ to the regular semisimple locus is a smooth morphism if $\mf{m}\supset\mf{b}$. 

It is well known \cite{Spr76} that the cohomology of a Springer fiber admits an action of the Weyl group $\msf{W}$ of $\msf{G}$, although $\msf{W}$ does not naturally act on the Springer fiber itself in general. Similarly, Tymoczko \cite{Tym08} showed that $\msf{W}$ acts on the cohomology of regular semisimple Hessenberg varieties (of type $A$) using the theory of moment graphs developed by Goresky, Kottwitz, and MacPherson \cite{GKM98}. According to \cite{BC18,BC24}, this action can be identified with the monodromy action arising from the local system obtained by the (derived) pushforward of the constant sheaf via $\mu^{\mr{rs}}_{\mf{m}}$. 

In type $A$, there is a bijection $\Gamma\mapsto\mf{n}_{\Gamma}$ between the set of unit interval graphs and the set of $\msf{B}$-stable subspaces of $\mf{n}$ given by
\begin{align*}
	\mf{n}_{\Gamma}=\Span\left\{ E_{ij}\mid i<j,\{i,j\}\notin\Gamma\right\},
\end{align*}
where $E_{ij}$ is the matrix unit, and we take $\msf{B}$ to be the subgroup of upper triangular matrices. The orthogonal complement $\mf{n}_{\Gamma}^{\perp}$ (with respect to the trace form) is given by 
\begin{align*}
	\mf{n}_{\Gamma}^{\perp}=\mf{b}+\Span\left\{ E_{ji}\mid i<j,\{i,j\}\in\Gamma\right\}.
\end{align*}
Let $D$ be a regular semisimple element, e.g., a diagonal matrix with distinct diagonal entries. We denote by $\msf{Hess}_{\Gamma,D}$ the regular semisimple Hessenberg variety associated with the subspace $\mf{m}=\mf{n}^{\perp}_{\Gamma}$ and the element $D$. In terms of the Hessenberg function $\msf{h}$ corresponding to $\Gamma$, $\msf{Hess}_{\Gamma,D}$ can be explicitly described in terms of flags of subspaces of $\bb{C}^n$ by
\begin{align*}
	\msf{Hess}_{\Gamma,D}\cong\left\{0\subset V_1\subset V_2\subset\cdots\subset V_{n}=\bb{C}^n\mid \dim V_i=i, DV_i\subset V_{\msf{h}(i)}\right\}.
\end{align*}

In \cite{SW12,SW16}, Shareshian and Wachs conjectured that 
\begin{align*}
	\omega\mb{X}_{\Gamma}(x;q)=\sum_{i\geq0}q^{i}\,\msf{ch}\left(H^{2i}(\msf{Hess}_{\Gamma,D})\right).
\end{align*}
Here, the cohomology $H^{\ast}(\msf{Hess}_{\Gamma,D})$ of $\msf{Hess}_{\Gamma,D}$ is equipped with Tymoczko's action of $\mf{S}_n$, and $\msf{ch}$ is the Frobenius characteristic of representations of $\mf{S}_n$. This conjecture was proved independently by Brosnan and Chow \cite{BC18} and Guay-Paquet \cite{GP16}, and has led to many studies on the geometry of Hessenberg varieties, such as \cite{AHMMS20,AN23,AN25b,AMS22,CHL23a,CHL23b,HP19,HP22,HPT21,HMS24,KL24a,KL24b,Les25,MS24,PS25}.

\subsubsection{Characters of Hecke algebras}

The \emph{Hecke algebra} $\mca{H}_{n}$ of type $A_{n-1}$ is a deformation of the group ring of $\mf{S}_{n}$. It is an algebra over $\bb{Z}[q^{1/2},q^{-1/2}]$ generated by $T_1,\ldots,T_{n-1}$ subject to the relations:
\begin{itemize}
	\item $(T_i-q^{1/2})(T_i+q^{-1/2})=0$ for $i=1,\ldots,n-1$,
	\item $T_{i}T_{i+1}T_{i}=T_{i+1}T_{i}T_{i+1}$ for $i=1,\ldots,n-2$, and 
	\item $T_iT_j=T_jT_i$ for $|i-j|\geq2$.
\end{itemize}
This algebra has a \emph{standard basis} $\{T_w\}_{w\in\mf{S}_n}$ over $\bb{Z}[q^{1/2},q^{-1/2}]$ indexed by $w\in\mf{S}_n$, defined by $T_w=T_{i_1}T_{i_2}\cdots T_{i_l}$ for a \emph{reduced} (i.e., shortest) expression $w=s_{i_1}s_{i_2}\cdots s_{i_l}$, where $s_i$ is the transposition of $i$ and $i+1$. The \emph{length} $l$ of such an expression is denoted by $\ell(w)$. 

There is a $\bb{Z}$-algebra involution $\beta$ on $\mca{H}_n$, called the \emph{bar involution}, given by
\begin{align*}
	\beta\left(\sum_{w\in\mf{S}_n}f_w(q)T_w\right)=\sum_{w\in\mf{S}_n}f_w(q^{-1})T_{w^{-1}}^{-1}
\end{align*}
for $f_w(q)\in\bb{Z}[q^{1/2},q^{-1/2}]$. In \cite{KL79}, Kazhdan and Lusztig introduced a remarkable basis $\{C_w\}_{w\in\mf{S}_n}$ of $\mca{H}_n$ called the \emph{Kazhdan--Lusztig basis}. This basis encodes significant information about the representation theory of semisimple Lie algebras, as established by the \emph{Kazhdan--Lusztig conjecture} (proved by Beilinson--Bernstein \cite{BB81} and Brylinski--Kashiwara \cite{BK81}). It is uniquely characterized by 
 \begin{itemize}
 	\item (\emph{bar invariance}) \, $\beta\left(C_w\right)=C_w$, and 
 	\item (\emph{triangularity}) \, $C_w\in T_w+\sum_{y<w}q^{-1/2}\bb{Z}[q^{-1/2}]\cdot T_y$.
 \end{itemize}
Here, the partial order $<$ in $\mf{S}_n$ is the \emph{Bruhat order}, i.e., $y\leq w$ if and only if an reduced expression for $y$ is a subexpression of a reduced expression for $w$. The theory of \emph{canonical bases}, characterized by analogues of the bar invariance and the triangularity, plays a fundamental role in the representation theory of quantum groups and various other algebras, and has been the subject of extensive research.

A connection between the theory of immanants and the Kazhdan--Lusztig basis was first observed by Haiman \cite{Hai93}. By Tits' deformation theorem (see for example \cite{GP00}), the generic representation theory of $\mca{H}_n$ is the same as the representation theory of $\mf{S}_n$. In particular, there are characters $\chi^{\lambda}_q$ (resp. $\phi^{\lambda}_q$) of $\mca{H}_n$ corresponding to the characters $\chi^{\lambda}$ (resp. $\phi^{\lambda}$) of $\mf{S}_n$. By definition, the immanant $\msf{Imm}_{\chi}H_{\mu/\nu}(x)$ with respect to a character $\chi\in\msf{R}_n$ can be written as $\chi(I_{\mu/\nu})$ for
\begin{align*}
	I_{\mu/\nu}\coloneqq\sum_{w\in\mf{S}_n}w\otimes h_{\mu+\delta-w(\nu+\delta)}(x)\in\bb{Z}[\mf{S}_n]\otimes\Lambda.
\end{align*}

Haiman \cite{Hai93} proved that $I_{\mu/\nu}$ is a nonnegative linear combination of $\left.C_w\right|_{q=1}\otimes s_{\theta}(x)$ for $w\in\mf{S}_n$ and $\theta\vdash N$, by interpreting each coefficient as the multiplicity of an irreducible $\mf{sl}_n$-representation in the tensor product of two irreducible representations, using the Kazhdan--Lusztig conjecture. Haiman also showed that $\chi^{\lambda}_q\left(q^{\ell(w)/2}C_w\right)\in\bb{Z}_{\geq0}[q]$ by combining the realization of the irreducible representations of $\mca{H}_n$ as cell representations (a feature specific to type $A$) \cite{KL79} with the positivity of the structure constants of the Kazhdan--Lusztig basis \cite{Spr82}. This in particular implies the $s$-positivity of $\msf{Imm}_{\chi^{\lambda}}H_{\mu/\nu}(x)$. Moreover, Haiman conjectured that $\phi^{\lambda}_q\left(q^{\ell(w)/2}C_w\right)\in\bb{Z}_{\geq0}[q]$.

In \cite{CHSS16}, Clearman, Hyatt, Shelton, and Skandera established the following formula for the chromatic quasisymmetric function of a unit interval graph:
\begin{align*}
	\mb{X}_{\Gamma}(x;q)=\sum_{\lambda\vdash n}\chi^{\lambda}_q\left(q^{\ell(w)/2}C_w\right)s_{\lambda^t}(x)=\sum_{\lambda\vdash n}\phi^{\lambda}_{q}\left(q^{\ell(w)/2}C_w\right)e_{\lambda}(x),
\end{align*}
where $w\in\mf{S}_n$ is the 312-avoiding permutation corresponding to $\Gamma$. In particular, Haiman's conjecture generalizes Conjecture~\ref{Conj_SW}. For recent progress in this direction, see for example \cite{AN24,AN25a,AN25b,Ska21}. For a similar story in type $BC$, see \cite{Ska25}.

\subsubsection{Other connections}

The interplay with combinatorics, geometry, and representation theory discussed above, along with further links to Macdonald polynomials \cite{HW20}, the KP hierarchy \cite{CKL20}, and characters of finite general linear groups \cite{Gag24}, to name a few, highlights the richness of the subject. We also mention that numerous variants of the chromatic (quasi)symmetric functions have been studied, for example, in \cite{CPS23,GS01,HOM25,Mar25,SY18,Sta98}. 

\subsection{Modular law}\label{Sec:Modular_Law}

In the study of chromatic symmetric functions, a certain linear relation between them, called the \textit{modular law}, plays a crucial role. This relation first appeared in the work of Guay-Paquet \cite[Proposition~3.1]{GP13} in the setting of $(\bf 3+1)$-free posets, and was generalized by Orellana and Scott \cite[Theorem~3.1]{OS14} to the so-called \textit{triple deletion property}. 

\subsubsection{Modular triples} 

A $q$-analogue of the modular law for unit interval graphs was first explicitly formulated by Lee \cite[Theorem~3.5]{Lee21} in the setting of unicellular LLT polynomials. Huh, Nam, and Yoo \cite[Theorem~3.1]{HNY20} reformulated this relation in the language of chromatic quasisymmetric functions using (reverse) area sequences. A dual version of Lee's relation was explicitly given by Alexandersson \cite[Proposition~23]{Ale21}, also within the framework of unicellular LLT polynomials. Let us state these relations in terms of conjugate Hessenberg functions. 

\begin{dfn}\label{Def:Modular_Triple}
A triple $(\msf{e},\msf{e}',\msf{e}'')$ in $\bb{E}_n$ is called a \emph{modular triple of type I} if there exists $i\in[n]$ such that $\msf{e}(i+1)=\msf{e}(i)+1$, 
\begin{align*}
	\msf{e}'(j)=\begin{cases}
		\msf{e}(i) & \text{if } j=i+1,\\
		\msf{e}(j) & \text{otherwise},
	\end{cases}\qquad\msf{e}''(j)=\begin{cases}
		\msf{e}(i+1) & \text{if } j=i,\\
		\msf{e}(j) & \text{otherwise},
	\end{cases}
\end{align*}
and such that $\msf{e}^{-1}(i)=\emptyset$.

A triple $(\msf{e},\msf{e}',\msf{e}'')$ in $\bb{E}_n$ is called a \emph{modular triple of type II} if there exists $i\in[n]$ such that 
\begin{align*}
	\msf{e}'(j)=\begin{cases}
		\msf{e}(i)-1 & \text{if } j=i,\\
		\msf{e}(j) & \text{otherwise},
	\end{cases}\quad\msf{e}''(j)=\begin{cases}
		\msf{e}(i)+1 & \text{if } j=i,\\
		\msf{e}(j) & \text{otherwise},
	\end{cases}
\end{align*} 
and such that $\msf{e}(\msf{e}(i))=\msf{e}(\msf{e}(i)+1)$. 
\end{dfn}

\begin{figure}
\begin{tikzpicture}
\node (type I) at (0,0) {
\begin{tikzpicture}[scale=0.8]
    \draw[step=1, gray!30] (0,0) grid (7,7);
    \draw[gray] (0,0) -- (7,7);
    \draw[very thick] 
        (0,0) -- (0,1) -- (0,2) -- (1,2) -- (1,3) -- (1,4) -- (2,4) -- (2,5) -- (3,5) -- (3,6) -- (4,6) -- (5,6)  -- (5,7) -- (6,7) -- (7,7);
    \draw[very thick, red] (1,4) -- (1,5) -- (2,5);
    \draw[very thick, blue] (1,3) -- (2,3) -- (2,4);
    \draw[very thick, dotted, ->, green] (2,4) -- (4,4) -- (4,6);
    \draw[ultra thick,green] (3.5,6) -- (4.5,6);
    \draw (-0.5,4) node {$i$};
    \draw (4,-0.5) node {$i$};
	\node at (3.5,-1.5) {Type I};
\end{tikzpicture}};
\node (type II) at (8,0) {
\begin{tikzpicture}[scale=0.8]
    \draw[step=1, gray!30] (0,0) grid (7,7);
    \draw[gray] (0,0) -- (7,7);
    \draw[very thick] 
        (0,0) -- (0,1) -- (0,2) -- (1,2) -- (1,3) -- (1,4) -- (2,4) -- (2,5) -- (3,5) -- (3,6) -- (4,6) -- (5,6) -- (5,7) -- (6,7) -- (7,7); 
    \draw[very thick, red] (2,5) -- (2,6) -- (3,6);
    \draw[very thick, blue] (3,5) -- (4,5) -- (4,6);
    \draw[very thick, dotted,->,green] (3,5) -- (3,3) -- (1,3);
    \draw[ultra thick,green] (1,2.5) -- (1,3.5);
    \draw (-0.5,6) node {$i$};
    \draw (3,-0.5) node {$\msf{e}(i)$};
    \node at (3.5,-1.5) {Type II};
\end{tikzpicture}};
\end{tikzpicture}
\caption{Modular triples $(\msf{e},\msf{e}',\msf{e}'')$ of type I and II. In both cases, $\msf{e}=(0,0,1,1,2,3,5)$ corresponds to the black path, $\msf{e}'$ corresponds to the red path, and $\msf{e}''$ corresponds to the blue path. The last conditions in Definition~\ref{Def:Modular_Triple} mean that the green dotted path intersects the Dyck path corresponding to $\msf{e}$ vertically.}\label{Fig_modular_triple}	
\end{figure}

An example illustrating these definitions is given in Figure~\ref{Fig_modular_triple}.

\begin{remark}\label{Rem:Modular_Triple}
If $(\msf{e},\msf{e}',\msf{e}'')$ is a modular triple of type I, then we have $i<n$ and $\msf{e}(i)+1=\msf{e}''(i)<i$. If $(\msf{e},\msf{e}',\msf{e}'')$ is a modular triple of type II, then we have $i>1$, $\msf{e}(i)+1=\msf{e}''(i)<i$, and $\msf{e}(i-1)<\msf{e}(i)<\msf{e}(i+1)$. 
\end{remark}

\begin{thm}[\cite{Ale21,Lee21}]\label{Thm:cqf_modular}
For any modular triple $(\msf{e},\msf{e}',\msf{e}'')$ of type I or II, we have
\begin{align*}
	\mb{X}_{\Gamma_{\msf{e}}}(x;q)=\frac{\mb{X}_{\Gamma_{\msf{e}'}}(x;q)+q\mb{X}_{\Gamma_{\msf{e}''}}(x;q)}{1+q}.
\end{align*}
\end{thm}

\begin{dfn}
A function $f:\bb{E}_n\rightarrow\bb{Q}(q)$ is said to satisfy the \emph{modular law} if 
\begin{align*}
	f(\msf{e})=\frac{f(\msf{e}')+qf(\msf{e}'')}{1+q}
\end{align*}
for any modular triple $(\msf{e},\msf{e}',\msf{e}'')$ of type I or II.
\end{dfn}

By Proposition~\ref{Prop_basic_properties} and Theorem~\ref{Thm:cqf_modular}, we obtain the following:

\begin{corollary}\label{Cor:c_lambda}
For any $\lambda\vdash n$, let $c_{\lambda}:\bb{E}_n\rightarrow\bb{Z}[q]$ be the function determined by
\begin{align*}
	\mb{X}_{\Gamma_{\msf{e}}}(x;q)=\sum_{\lambda\vdash n}c_{\lambda}(\msf{e};q)e_{\lambda}(x).
\end{align*}
Then the function $c_{\lambda}$ satisfies the modular law and 
\begin{align*}
	c_{\lambda}(\msf{e}_{\alpha};q)=\begin{cases}
		\prod_{i}[\lambda_i]_q! & \text{if }\lambda=(\lambda_1,\lambda_2,\ldots)\text{ is the weakly decreasing rearrangement of }\alpha,\\
		0 & \text{otherwise.}
	\end{cases}
\end{align*}
\end{corollary}

\subsubsection{Abreu--Nigro's characterization}

A crucial ingredient in our proof of the Stanley--Stembridge conjecture is a result by Abreu and Nigro \cite{AN21}, which characterizes chromatic quasisymmetric functions for unit interval graphs in terms of type I and II modular relations, along with the base cases for disjoint ordered unions of complete graphs. 

In order to state their results, let us introduce a family of conjugate Hessenberg functions corresponding to disjoint ordered unions of complete graphs.

\begin{dfn}\label{Def:Parabolic}
For a composition $\alpha=(\alpha_1,\alpha_2,\ldots,\alpha_l)$ of $n$, we define the associated \emph{parabolic conjugate Hessenberg function} $\msf{e}_{\alpha}\in\bb{E}_n$ by 
\begin{align*}
	\msf{e}_{\alpha}(i)=\sum_{j=1}^{k-1}\alpha_j
\end{align*}
for any $i\in I^{\alpha}_{k}$, where the intervals $I^{\alpha}_1,\ldots,I^{\alpha}_l$ partition the set $[n]$ and are defined by
\begin{align*}
	I^{\alpha}_k\coloneqq\biggl\{ i\in[n] \bigg| \sum_{j=1}^{k-1}\alpha_j < i \leq \sum_{j=1}^{k}\alpha_j \biggr\}.
\end{align*}
\end{dfn}

For example, $\msf{e}_{(n)}=(0,0,\ldots,0)$, and $\Gamma_{\msf{e}_{(n)}}=K_n$ is the complete graph $K_n$ on $n$ vertices. As an unlabeled graph, $\Gamma_{\msf{e}_{\alpha}}$ is isomorphic to the disjoint union $\bigsqcup_{i}K_{\alpha_i}$. The name comes from the fact that the corresponding $\msf{B}$-stable subspace $\mf{n}_{\Gamma}^{\perp}\subset\mf{gl}_n$ (see Section~\ref{Sec:Hess}) is a parabolic subalgebra of $\mf{gl}_n$.

\begin{thm}[{\cite[Theorem~1.2]{AN21}}]\label{Thm:AbrNig}
Let $f:\bb{E}_n\rightarrow\bb{Q}(q)$ be a function satisfying the modular law. Then $f$ is determined by its values $f(\msf{e}_{\alpha})$ at the parabolic conjugate Hessenberg functions $\msf{e}_{\alpha}$ for all $\alpha\models n$. Moreover, we have $f(\msf{e}_{\alpha})=f(\msf{e}_{\lambda})$ for the weakly decreasing rearrangement $\lambda$ of $\alpha$.
\end{thm}

Combining this with Corollary~\ref{Cor:c_lambda}, we obtain the following description of a function satisfying the modular law.

\begin{corollary}[{\cite[Corollary~3.2]{AN21}}]\label{Cor:f_to_c}
For a function $f:\bb{E}_n\rightarrow\bb{Q}(q)$ satisfying the modular law, we have
\begin{align*}
	f(\msf{e})=\sum_{\lambda\vdash n}f(\msf{e}_{\lambda})\frac{c_{\lambda}(\msf{e};q)}{\prod_{i}[\lambda]_q!}
\end{align*}
for any $\msf{e}\in\bb{E}_n$.
\end{corollary}

As an easy corollary of Theorem~\ref{Thm:AbrNig}, we obtain the following relation.

\begin{corollary}\label{Cor:Prob_Rel}
For any $\msf{e}\in\bb{E}_n$, we have
\begin{align*}
	\sum_{\lambda\vdash n}\frac{c_{\lambda}(\msf{e};q)}{\prod_{i}[\lambda_i]_q!}=1.
\end{align*}
\end{corollary}

\begin{proof}
This follows from Corollary~\ref{Cor:f_to_c} applied to the constant function $f=1$, which satisfies the modular law.
\end{proof}

\subsection{Kato's formula}

The starting point of this paper is Kato's formula \cite{Kat24} for the chromatic symmetric function of a unit interval graph $\Gamma$. This work provides a new geometric interpretation of $\mb{X}_{\Gamma}(x;q)$, and also establishes a new algebraic formula for $\mb{X}_{\Gamma}(x;q=1)$ in terms of the affine Weyl group of type $A$. 

\subsubsection{Geometric part}

Recall from Section~\ref{Sec:Hess} that $\Gamma$ corresponds to a $\msf{B}$-stable subspace $\mf{n}_{\Gamma}\subset\mf{n}$ for $\msf{G}=\GL_n(\bb{C})$. Kato \cite{Kat23} constructed a smooth compactification $\mf{X}_{\Gamma}$ of $\msf{G}\times^{\msf{B}}\mf{n}_{\Gamma}$ equipped with an action of $\msf{G}(\msc{O})$ (where $\msc{O}\coloneqq\bb{C}[\![z]\!]$ is the ring of formal power series over $\bb{C}$), and a $\msf{G}(\msc{O})$-equivariant proper morphism $\msf{m}_{\Gamma}:\mf{X}_{\Gamma}\rightarrow\msf{Gr}\coloneqq\msf{G}(\msc{K})/\msf{G}(\msc{O})$ into the \emph{affine Grassmannian} $\msf{Gr}$ of type $A$ (where $\msc{K}\coloneqq\bb{C}(\!(z)\!)$ is the field of Laurent series over $\bb{C}$). 

By the decomposition theorem \cite{BBD81}, the (derived) pushforward $\msf{m}_{\Gamma\ast}\underline{\bb{C}}_{\mf{X}_{\Gamma}}$ of the constant sheaf on $\mf{X}_{\Gamma}$ by $\msf{m}_{\Gamma}$ decomposes into a direct sum of shifted $\msf{G}(\msc{O})$-equivariant irreducible perverse sheaves on $\msf{Gr}$, i.e., there exist non-negative integers $m^{\Gamma}_{\lambda,i}$ for each dominant cocharacter $\lambda$ of $\msf{G}$ and $i\in\bb{Z}$ such that
\begin{align*}
	\msf{m}_{\Gamma\ast}\underline{\bb{C}}_{\mf{X}_{\Gamma}}\cong\bigoplus_{\substack{\lambda:\bb{G}_{\mr{m}}\rightarrow\msf{G}\\\lambda:\text{ dominant}}}\bigoplus_{i\in\bb{Z}}\msf{IC}_{\lambda}[-i]^{m^{\Gamma}_{\lambda,i}},
\end{align*}
where $\msf{IC}_{\lambda}$ is the $\msf{G}(\msc{O})$-equivariant irreducible perverse sheaf on $\msf{Gr}$ obtained by the intermediate extension of the shifted constant sheaf on the $\msf{G}(\msc{O})$-orbit containing the image of $z\in\bb{G}_{\mr{m}}(\msc{K})$ under $\bb{G}_{\mr{m}}(\msc{K})\xrightarrow{\lambda}\msf{G}(\msc{K})\rightarrow\msf{Gr}$. 

By the geometric Satake correspondence \cite{MV07}, the category of $\msf{G}(\msc{O})$-equivariant perverse sheaves on $\msf{Gr}$ is equivalent to the category of representations of $\GL_n(\bb{C})$. The perverse sheaf $\msf{IC}_{\lambda}$ corresponds to the irreducible highest weight representation of $\GL_n(\bb{C})$ with highest weight $\lambda$, whose character is given by the Schur polynomial $s_{\lambda}(x_1,\ldots,x_n)$. Kato \cite{Kat24} obtained the following interpretation of (the $n$-th truncation of) $\mb{X}_{\Gamma}(x;q)$ as the character of the graded $\GL_n(\bb{C})$-representation corresponding to $\msf{m}_{\Gamma\ast}\underline{\bb{C}}_{\mf{X}_{\Gamma}}$:
\begin{align}\label{Eqn:Kato_Geom}
	\mb{X}_{\Gamma}(x_1,\ldots,x_n;q)=q^{-d_{\Gamma}}\sum_{\lambda:\text{ dominant}}\sum_{i\in\bb{Z}}m^{\Gamma}_{\lambda,i}q^{i/2}s_{\lambda}(x_1,\ldots,x_n),
\end{align} 
where $d_{\Gamma}\in\bb{Z}_{\geq0}$ is an explicitly determined constant. The proof in \cite{Kat24} is based on Abreu--Nigro's characterization \cite{AN21}. Kamnitzer \cite{Kam24} gave another direct proof of this equality.

\subsubsection{Algebraic part}

Let $\msf{e}\in\bb{E}_n$ be the conjugate Hessenberg function corresponding to $\Gamma$. Another interesting feature of Kato's variety $\mf{X}_{\Gamma}$ is that it is constructed as an iterated projective bundle. This fact, together with (\ref{Eqn:Kato_Geom}), implies (via a localization calculation \cite{Kat24}) that $\mb{X}_{\Gamma}(x_1,\ldots,x_n;q=1)$ is obtained as an iterated application of operators
\begin{align*}
	\msf{S}_{i,\msf{e}(i)}\coloneqq1+s_i+s_{i+1}s_i+\cdots+s_{n-1}\cdots s_i+s_0s_{n-1}\cdots s_i+\cdots+s_{\msf{e}(i)-1}\cdots s_0s_{n-1}\cdots s_i
\end{align*}
in the group ring of the affine Weyl group $\widetilde{\mf{S}}_n=\langle s_i\mid i\in\bb{Z}/n\bb{Z}\rangle$ of type $A_{n-1}$ to $x_1\cdots x_n\in\bb{Q}[x_1^{\pm1},\ldots,x_n^{\pm1}]$, i.e.,
\begin{align}\label{Eqn:Kato_Alg}
	\mb{X}_{\Gamma}(x_1,\ldots,x_n)=\msf{S}_{1,\msf{e}(1)}\msf{S}_{2,\msf{e}(2)}\cdots \msf{S}_{n,\msf{e}(n)}(x_1x_2\cdots x_n).
\end{align}
Here, the action of $s_i\in \widetilde{\mf{S}}_n$ on $\bb{Q}[x_1^{\pm1},\ldots,x_n^{\pm1}]$ is given by the transposition of $x_i$ and $x_{i+1}$ for $1\leq i\leq n-1$ and 
\begin{align}\label{Eqn:s0}
	s_0\left(f(x_1,\ldots,x_n)\right)=\frac{x_1}{x_n}f(x_n,x_2,\ldots,x_{n-1},x_1).
\end{align}

It is straightforward \cite{Hik25b} to introduce the parameter $q$ into (\ref{Eqn:Kato_Alg}) by considering the \emph{affine} Hecke algebra of type $A_{n-1}$. It is also natural to introduce another parameter (as done in \cite{Kat24}) into (\ref{Eqn:s0}) from the point of view of affine Lie algebras \cite{Kac90}. This leads us to define a $(q,t)$-analogue of the chromatic symmetric functions in \cite{Hik25b}\footnote{The role of $q$ here is played by $t$ in \cite{Hik25b}.}. 

\begin{remark}
We arrived at the statement of Corollary~\ref{Cor:Prob_Rel} while looking for an application of our $(q,t)$-chromatic symmetric functions in \cite{Hik25b}, where it is stated in a slightly different form. They are related by $q\mapsto q^{-1}$ via the property (vi) in Proposition~\ref{Prop_basic_properties}.
\end{remark}

\begin{remark}
The author (jointly with Syu Kato) previously attempted to approach the Stanley--Stembridge conjecture using (\ref{Eqn:Kato_Alg}). In the course of this investigation, the author was led to speculate that a formula with operators applied in the reverse direction would be more amenable to inductive arguments. Our probabilistic formula for $\mb{X}_{\Gamma}(x;q)$ takes precisely this form (see Definition~\ref{Def:f_lambda}). 
\end{remark}

\section{Reduction to Key Lemmas}\label{Sec:Reduction}

We continue to use the notation introduced in Section~\ref{Sec:Main}. We now proceed to the proof of our main theorem (Theorem~\ref{Thm:Main}). By Abreu--Nigro's characterization (Theorem~\ref{Thm:AbrNig}), it suffices to check that for any $\lambda\vdash n$, the function $\msf{e}\mapsto \Prob\{\msf{X}^{(\msf{e})}_n=\lambda\}$ satisfies the modular law and takes the correct values at parabolic conjugate Hessenberg functions.

The main task in the proof is to verify the modular law. Recall the definition of a modular triple $(\msf{e},\msf{e}',\msf{e}'')$ in Definition~\ref{Def:Modular_Triple}. While $\msf{e}'$ and $\msf{e}''$ are obtained from $\msf{e}$ by local modification at $i$ or $i+1$, the definition involves non-local conditions, namely $\msf{e}^{-1}(i)=\emptyset$ or $\msf{e}(\msf{e}(i))=\msf{e}(\msf{e}(i)+1)$. 

In this section, we first reformulate our main theorem in linear algebraic terms and check the values at parabolic conjugate Hessenberg functions. We then state several key lemmas concerning the properties of operators derived from the transition probabilities of the time-inhomogeneous Markov chain $\{\widetilde{\msf{X}}^{(\msf{e})}_n\}$. These properties arise from the study of local modifications in a modular triple and should be regarded as an \emph{operator-level modular law}. We reduce the proof of our main theorem to these key lemmas by demonstrating how this operator-level modular law propagates to the actual modular law under the non-local conditions satisfied by modular triples.

\subsection{Reformulation}\label{Sec:Reformulation}

Since the stochastic process $\{\widetilde{\msf{X}}^{(\msf{e})}_{n}\}_{n\in\bb{Z}_{\geq0}}$ in Definition~\ref{Def:Stochastic_Process} is a Markov chain, the probability $\Prob\bigl\{\widetilde{\msf{X}}^{(\msf{e})}_{n}=T\bigr\}$ for any $T\in\msf{SYT}_n$ can be computed by taking the product of the operators determined by the transition probabilities. 

\begin{dfn}
Let $\msf{V}_n\coloneqq\Span_{\bb{Q}(q)}\left\{T\mid T\in\msf{SYT}_n\right\}$ be the $\bb{Q}(q)$-vector space formally generated by the set of standard Young tableaux of size $n$. For each $r\in\bb{Z}_{\geq0}$, we define a $\bb{Q}(q)$-linear map $\Omega_r:\msf{V}_{n-1}\rightarrow\msf{V}_n$ by
\begin{align*}
	\Omega_r(T)&\coloneqq\sum_{T'\in\msf{SYT}_n}\msf{P}^{(r)}_n(T;T')\cdot T'\\
	&=\sum_{c\in W\left(\boldsymbol{\delta}^{(r)}(T)\right)}\varphi_c\bigl(\boldsymbol{\delta}^{(r)}(T);q\bigr)\cdot \mf{f}_c(T)
\end{align*}
for any $T\in\msf{SYT}_{n-1}$. 
\end{dfn}

\begin{dfn}\label{Def:f_lambda}
For each $\lambda\vdash n$, we define a function $\msf{f}_{\lambda}:\bb{E}_n\rightarrow\bb{Q}(q)$ by
\begin{align*}
	\msf{f}_{\lambda}(\msf{e})\coloneqq\pi_{\lambda}\Omega_{\msf{e}(n)}\Omega_{\msf{e}(n-1)}\cdots\Omega_{\msf{e}(1)}(\emptyset),
\end{align*}
where $\pi_{\lambda}:\msf{V}_n\rightarrow\bb{Q}(q)$ is the linear map defined by 
\begin{align*}
	\pi_{\lambda}(T)=\begin{cases}
		1 & \text{if }T\in\msf{SYT}(\lambda),\\
		0 & \text{otherwise},
	\end{cases}
\end{align*}
for any $T\in\msf{SYT}_n$.
\end{dfn}

Using these operators, we obtain
\begin{align*}
	\Omega_{\msf{e}(n)}\Omega_{\msf{e}(n-1)}\cdots\Omega_{\msf{e}(1)}(\emptyset)=\sum_{T\in\msf{SYT}_n}\Prob\bigl\{\widetilde{\msf{X}}^{(\msf{e})}_{n}=T\bigr\}\cdot T
\end{align*}
for any $\msf{e}\in\bb{E}_n$, and hence 
\begin{align*}
	\Prob\left\{\msf{X}^{(\msf{e})}_n=\lambda\right\}=\msf{f}_{\lambda}(\msf{e})
\end{align*}
for any $\lambda\vdash n$. 

Recall the notion of the parabolic conjugate Hessenberg function $\msf{e}_{\alpha}\in\bb{E}_n$ for a composition $\alpha\models n$ from Definition~\ref{Def:Parabolic}. By Corollary~\ref{Cor:f_to_c}, Theorem~\ref{Thm:Main} is reduced to the following properties of the function $\msf{f}_{\lambda}$.

\begin{thm}\label{Thm:Ref}
For each $\lambda\vdash n$, the function $\msf{f}_{\lambda}:\bb{E}_n\rightarrow\bb{Q}(q)$ satisfies the modular law and 
\begin{align*}
	\msf{f}_{\lambda}(\msf{e}_{\alpha})=\begin{cases}
		1 & \text{if }\lambda\text{ is the weakly decreasing rearrangement of }\alpha,\\
		0 & \text{otherwise},
	\end{cases}
\end{align*}
for any composition $\alpha\models n$.
\end{thm}

\begin{proof}
The proof is organized as follows. We establish the second assertion concerning the values at $\msf{e}_{\alpha}$ in Proposition~\ref{Prop:Parabolic} in Section~\ref{Sec:Parabolic}. In Section~\ref{Sec:Operator}, we verify the type I modular law in Proposition~\ref{Prop:Modular_I} and the type II modular law in Proposition~\ref{Prop:Modular_II}, assuming the validity of Lemmas~\ref{Lem:Key1}, \ref{Lem:Key2}, and \ref{Lem:Key3}. Finally, we prove these key lemmas in Section~\ref{Sec:Key}.
\end{proof}

\begin{remark}
Given that $\msf{f}_{\lambda}$ satisfies the modular law, the value of $\msf{f}_{\lambda}$ at $\msf{e}_{\alpha}$ for $\alpha\models n$ coincides with the value at $\msf{e}_{\mu}$ for the weakly decreasing rearrangement $\mu\vdash n$ of $\alpha$, by Theorem~\ref{Thm:AbrNig}. Hence it suffices to consider the case where $\alpha\vdash n$ in Theorem~\ref{Thm:Ref}.
\end{remark}

\subsection{Parabolic cases}\label{Sec:Parabolic}

In this section, we prove the second statement of Theorem~\ref{Thm:Ref}. 

We construct the standard Young tableau $T_{\alpha}\in\msf{SYT}(\mu)$ for each composition $\alpha\models n$, where $\mu\vdash n$ is its weakly decreasing rearrangement. First, we identify the composition $\alpha$ with its diagram $\{(i,j)\in\bb{Z}^2_{>0}\mid 1\leq i\leq \alpha_j\}$ (in French notation). We fill the boxes of this diagram with the numbers $1,2,\ldots,n$ in increasing order, proceeding from bottom to top and from left to right within each row. We then obtain the standard Young tableau $T_{\alpha}$ of shape $\mu$ by letting the boxes fall under gravity. See Figure~\ref{Fig:Gravity} for an example.

\begin{figure}
\begin{tikzpicture}
\node at (0,0) {
\begin{ytableau}[] 
13 & 14 & 15\\
8 & 9 & 10 & 11 & 12\\
7\\
3 & 4 & 5 & 6\\
1 & 2 
\end{ytableau}};
\node at (4.75,0) {$T_{\alpha}=$};
\node at (7,0) {
\begin{ytableau}[] 
13\\
8 & 14\\
7 & 9 & 15\\
3 & 4 & 10 & 11\\
1 & 2 & 5 & 6 & 12
\end{ytableau}};
\draw[->,thick] (2.3,0) -- (4,0);
\end{tikzpicture}
\caption{$T_{\alpha}$ for $\alpha=(2,4,1,5,3)$.}\label{Fig:Gravity}
\end{figure}

\begin{prop}\label{Prop:Parabolic}
For any $\alpha\models n$, we have
\begin{align*}
	\Omega_{\msf{e}_{\alpha}(n)}\Omega_{\msf{e}_{\alpha}(n-1)}\cdots\Omega_{\msf{e}_{\alpha}(1)}(\emptyset)=T_{\alpha}.
\end{align*}
\end{prop}

\begin{proof}
Let $\alpha=(\alpha_1,\ldots,\alpha_l)$ be a composition. Recall the partition of the set $[n]=I^{\alpha}_1\sqcup\cdots\sqcup I^{\alpha}_{l}$ given in Definition~\ref{Def:Parabolic}. We divide the $n$ steps of adding boxes into $l$ stages according to this set partition. 

Since the definition of $\Omega_r$ depends only on the upper part of the standard Young tableau, the actual heights of the boxes added during the process are irrelevant, provided that they are eventually dropped down. Thus, in the $k$-th stage, we may assume that there is a floor at height $k-1$, and that boxes do not fall below this level.

In the first stage, we have $\msf{e}_{\alpha}(i)=0$ for $i\in I^{\alpha}_1=\{1,\ldots,\alpha_1\}$. At the $j$-th step of this stage, the boxes with labels greater than $\msf{e}_{\alpha}(j)=0$ occupy the first $j-1$ columns on the floor at height 0. Therefore, the corresponding Maya diagram is of the form $\delta=(1^{\infty},0^{\infty})$ with $W(\delta)=\{j\}$. Hence the operator $\Omega_{\msf{e}_{\alpha}(j)}$ adds a new box at the $j$-th column with probability $1$. At the end of the first stage, we obtain the unique standard Young tableau of shape $(\alpha_1)$. 

For the $k$-th stage, note that $\msf{e}_{\alpha}(i)=\max I^{\alpha}_{k-1}$ for any $i\in I^{\alpha}_k$. This implies that the boxes added in the previous stages (1 to $k-1$) are not colored red in stage $k$. Thus we can proceed analogously to the first stage: at the $j$-th step of stage $k$, the new box is added at the $j$-th column on the floor at height $k-1$ with probability 1.

At the end of stage $l$, the configuration of boxes corresponds to the skew diagram of $\alpha$ with labels $1,2,\ldots,n$ in increasing order from bottom to top and from left to right within each row. Finally, by removing the floors and letting the boxes fall, we obtain the standard Young tableau $T_{\alpha}$ with probability $1$. 
\end{proof}

\begin{corollary}\label{Cor:Parabolic}
For any $\lambda\vdash n$ and $\alpha\models n$, we have
\begin{align*}
	\msf{f}_{\lambda}(\msf{e}_{\alpha})=\begin{cases}
		1 & \text{if }\lambda\text{ is the weakly decreasing rearrangement of }\alpha,\\
		0 & \text{otherwise}.
	\end{cases}
\end{align*}
\end{corollary}

\begin{proof}
This follows immediately from Proposition~\ref{Prop:Parabolic} and 
\begin{align*}
	\pi_{\lambda}(T_{\alpha})=\begin{cases}
		1 & \text{if }\lambda\text{ is the weakly decreasing rearrangement of }\alpha,\\
		0 & \text{otherwise}.
		\end{cases}
\end{align*}
\end{proof}

\subsection{Auxiliary subspaces}\label{Sec:Auxiliary}

In order to clarify the structure of the argument, we introduce two auxiliary subspaces $\msf{K}_{m,n},\msf{M}_{m,n}\subset\msf{V}_n$ for $m\in[n-1]$ as follows.

\begin{dfn}
For $T\in\msf{SYT}_n$ and $m\in[n-1]$, let $T'$ be the tableau obtained from $T$ by swapping the boxes containing $m$ and $m+1$. If $T'$ is a standard Young tableau, then we define $\tau_m(T)\coloneqq T'$.
\end{dfn}

We will frequently use the obvious identity $\mf{f}_c\left(\tau_m(T)\right)=\tau_m\left(\mf{f}_c(T)\right)$ whenever both sides are well-defined.

\begin{dfn}\label{Def:mu}
For $T\in\msf{SYT}_n$ and $m\in[n-1]$, we define $\mu_m(T)\in\msf{V}_n$ by 
\begin{align*}
	\mu_m(T)\coloneqq\begin{cases}
		\frac{[a-b-1]_q}{[a-b]_q}\cdot T+\frac{[b-a-1]_q}{[b-a]_q}\cdot\tau_m(T) &\text{if }\msf{top}_T(a)=m \text{ and }\msf{top}_T(b)=m+1,\\
		0 &\text{otherwise}.
	\end{cases}
\end{align*}
Here, we recall that $\msf{top}_T(i)$ is the entry in the topmost box of the $i$-th column of $T$.
\end{dfn}

We note that if $\msf{top}_T(a)=m$, $\msf{top}_T(b)=m+1$, and $\tau_m(T)$ is not well-defined, then we have $b=a+1$. In such a case, we set $\mu_m(T)\coloneqq[2]_{q^{-1}}\cdot T$. 

\begin{dfn}\label{Def:Auxiliary}
For $m\in[n-1]$, we define two subspaces $\msf{K}_{m,n},\msf{M}_{m,n}\subset\msf{V}_n$ by
\begin{align*}
	\msf{K}_{m,n}&\coloneqq\Span_{\bb{Q}(q)}\left\{T-\tau_m(T)\mid T\in\msf{SYT}_n, \tau_m(T)\text{: well-defined}\right\},\\
	\msf{M}_{m,n}&\coloneqq\Span_{\bb{Q}(q)}\left\{\mu_m(T)\mid T\in\msf{SYT}_n\right\}.
\end{align*}
\end{dfn}

We obviously have $\msf{K}_{m,n}\subset\Ker\pi_{\lambda}$ for any $\lambda\vdash n$. In addition, we have the following basic properties.

\begin{lemma}\label{Lem:Basic}
Let $m\in[n-1]$.
\begin{enumerate}
	\item \label{item:i} If $r\neq m$, then we have $\Omega_r\left(\msf{K}_{m,n}\right)\subset\msf{K}_{m,n+1}$.
	\item \label{item:ii} If $r<m$, then we have $\Omega_r\left(\msf{M}_{m,n}\right)\subset\msf{M}_{m,n+1}$.
\end{enumerate}
\end{lemma}

\begin{proof}
If $r\neq m$, then we have $\boldsymbol{\delta}^{(r)}(T)=\boldsymbol{\delta}^{(r)}\left(\tau_m(T)\right)\eqqcolon\delta$ for any $T\in\msf{SYT}_n$. This implies
\begin{align*}
	\Omega_r(T)-\Omega_r\left(\tau_m(T)\right)&=\sum_{c\in W(\delta)}\varphi_c(\delta)\left(\mf{f}_c(T)-\mf{f}_c\left(\tau_m(T)\right)\right)\\
&=\sum_{c\in W(\delta)}\varphi_c(\delta)\left(\mf{f}_c(T)-\tau_m\left(\mf{f}_c(T)\right)\right)\in \msf{K}_{m,n+1}.
\end{align*} 
This proves the first assertion. 

For the second assertion, we assume $\msf{top}_T(a)=m$, $\msf{top}_T(b)=m+1$, and $r<m$. Then, since the entries $m$ and $m+1$ are greater than $r$, the corresponding columns $a$ and $b$ satisfy $\delta_a=\delta_b=1$. Thus we have $a,b\notin W(\delta)$. Consequently,  for any $c\in W(\delta)$, the operator $\mf{f}_c$ adds a box at a column different from $a$ and $b$, so $\mf{f}_c(T)$ also satisfies $\msf{top}_{\mf{f}_c(T)}(a)=m$ and $\msf{top}_{\mf{f}_c(T)}(b)=m+1$. This implies
\begin{align*}
	\Omega_r\left(\mu_m(T)\right)&=\frac{[a-b-1]_q}{[a-b]_q}\cdot \Omega_r(T)+\frac{[b-a-1]_q}{[b-a]_q}\cdot\Omega_r\left(\tau_m(T)\right)\\
	&=\sum_{c\in W(\delta)}\varphi_c(\delta)\left(\frac{[a-b-1]_q}{[a-b]_q}\cdot \mf{f}_c(T)+\frac{[b-a-1]_q}{[b-a]_q}\cdot\mf{f}_c\left(\tau_m(T)\right)\right)\\
	&=\sum_{c\in W(\delta)}\varphi_c(\delta)\mu_m\left(\mf{f}_c(T)\right)\in\msf{M}_{m,n+1}.
\end{align*}
This completes the proof of the second assertion.
\end{proof}

\subsection{Operator-level modular law}\label{Sec:Operator}

We now state three key lemmas, whose proofs are postponed to Section~\ref{Sec:Key}, and use them to verify the modular law for $\msf{f}_{\lambda}$. 

\begin{restatable}{lemma}{KeyLemmaOne}\label{Lem:Key1}
For any $0\leq r<n$, we have 
\begin{align*}
	\Omega_r\Omega_r(\msf{V}_{n-1})\subset\msf{M}_{n,n+1}.
\end{align*}
\end{restatable}

\begin{restatable}{lemma}{KeyLemmaTwo}\label{Lem:Key2}
For any $r\in[n-1]$, we have
\begin{align*}
	\left(\Omega_r-\frac{\Omega_{r-1}+q\Omega_{r+1}}{1+q}\right)(\msf{M}_{r,n})\subset\msf{K}_{r,n+1}.
\end{align*}
\end{restatable}

\begin{restatable}{lemma}{KeyLemmaThree}\label{Lem:Key3}
For any $0\leq r<n-1$, we have
\begin{align*}
	\left(\Omega_{r+1}\Omega_r-\frac{\Omega_r\Omega_r+q\Omega_{r+1}\Omega_{r+1}}{1+q}\right)(\msf{V}_{n-1})\subset\msf{K}_{n,n+1}.
\end{align*}
\end{restatable}

\begin{prop}\label{Prop:Modular_I}
For any $\lambda\vdash n$, the function $\msf{f}_{\lambda}:\bb{E}_n\rightarrow\bb{Q}(q)$ satisfies the modular law of type I.
\end{prop}

\begin{proof}
Let $(\msf{e},\msf{e}',\msf{e}'')$ be a modular triple of type I (see Definition~\ref{Def:Modular_Triple}). Let $i\in[n]$ be as in Definition~\ref{Def:Modular_Triple}. By Remark~\ref{Rem:Modular_Triple}, we have $\msf{e}(i)<i-1$ and hence we may apply Lemma~\ref{Lem:Key3} with $r=\msf{e}(i)$ and $n=i$ to obtain
\begin{align*}
\left(\Omega_{\msf{e}(i+1)}\Omega_{\msf{e}(i)}-\frac{\Omega_{\msf{e}'(i+1)}\Omega_{\msf{e}'(i)}+q\Omega_{\msf{e}''(i+1)}\Omega_{\msf{e}''(i)}}{1+q}\right)\Omega_{\msf{e}(i-1)}\cdots\Omega_{\msf{e}(1)}(\emptyset)\in\msf{K}_{i,i+1}.
\end{align*}
By the condition $\msf{e}^{-1}(i)=\emptyset$, we may apply Lemma~\ref{Lem:Basic}~\eqref{item:i} repeatedly to obtain
\begin{align*}
	\Omega_{\msf{e}(n)}\cdots\Omega_{\msf{e}(i+2)}\left(\Omega_{\msf{e}(i+1)}\Omega_{\msf{e}(i)}-\frac{\Omega_{\msf{e}'(i+1)}\Omega_{\msf{e}'(i)}+q\Omega_{\msf{e}''(i+1)}\Omega_{\msf{e}''(i)}}{1+q}\right)\Omega_{\msf{e}(i-1)}\cdots\Omega_{\msf{e}(1)}(\emptyset)\in\msf{K}_{i,n}.
\end{align*}
By applying $\pi_{\lambda}$, we obtain the modular relation 
\begin{align*}
	\msf{f}_{\lambda}(\msf{e})=\frac{\msf{f}_{\lambda}(\msf{e}')+q\msf{f}_{\lambda}(\msf{e}'')}{1+q},
\end{align*}
since $\pi_{\lambda}(\msf{K}_{i,n})=0$ for any $\lambda\vdash n$.
\end{proof}

\begin{prop}\label{Prop:Modular_II}
For any $\lambda\vdash n$, the function $\msf{f}_{\lambda}:\bb{E}_n\rightarrow\bb{Q}(q)$ satisfies the modular law of type II.
\end{prop}

\begin{proof}
Let $(\msf{e},\msf{e}',\msf{e}'')$ be a modular triple of type II, and let $i\in[n]$ be as in Definition~\ref{Def:Modular_Triple}. Since $\msf{e}(\msf{e}(i)+1)=\msf{e}(\msf{e}(i))$, we may apply Lemma~\ref{Lem:Key1} with $r=\msf{e}(\msf{e}(i))$ and $n=\msf{e}(i)$ to obtain
\begin{align*}
	\Omega_{\msf{e}(\msf{e}(i)+1)}\Omega_{\msf{e}(\msf{e}(i))}\cdots\Omega_{\msf{e}(1)}(\emptyset)\in\msf{M}_{\msf{e}(i),\msf{e}(i)+1}.
\end{align*}
By Remark~\ref{Rem:Modular_Triple}, we have $\msf{e}(i-1)<\msf{e}(i)<\msf{e}(i+1)$ and hence we may apply Lemma~\ref{Lem:Basic}~\eqref{item:ii} repeatedly to obtain
\begin{align*}
	\Omega_{\msf{e}(i-1)}\cdots\Omega_{\msf{e}(\msf{e}(i)+1)}\Omega_{\msf{e}(\msf{e}(i))}\cdots\Omega_{\msf{e}(1)}(\emptyset)\in\msf{M}_{\msf{e}(i),i-1}.
\end{align*}
Again by Remark~\ref{Rem:Modular_Triple}, we have $\msf{e}(i)<i-1$ and hence we may apply Lemma~\ref{Lem:Key2} for $r=\msf{e}(i)$ and $n=i-1$ to obtain
\begin{align*}
	\left(\Omega_{\msf{e}(i)}-\frac{\Omega_{\msf{e}'(i)}+q\Omega_{\msf{e}''(i)}}{1+q}\right)\Omega_{\msf{e}(i-1)}\cdots\Omega_{\msf{e}(1)}(\emptyset)\in\msf{K}_{\msf{e}(i),i}.
\end{align*}
By $\msf{e}(i+1)>\msf{e}(i)$, we may apply Lemma~\ref{Lem:Basic}~\eqref{item:i} repeatedly to obtain
\begin{align*}
	\Omega_{\msf{e}(n)}\cdots\Omega_{\msf{e}(i+1)}\left(\Omega_{\msf{e}(i)}-\frac{\Omega_{\msf{e}'(i)}+q\Omega_{\msf{e}''(i)}}{1+q}\right)\Omega_{\msf{e}(i-1)}\cdots\Omega_{\msf{e}(1)}(\emptyset)\in\msf{K}_{\msf{e}(i),n}.
\end{align*}
By applying $\pi_{\lambda}$, we obtain
\begin{align*}
	\msf{f}_{\lambda}(\msf{e})=\frac{\msf{f}_{\lambda}(\msf{e}')+q\msf{f}_{\lambda}(\msf{e}'')}{1+q}
\end{align*}
for any $\lambda\vdash n$.
\end{proof}

\section{Proof of the Key Lemmas}\label{Sec:Key}

In this section, we complete the proof of our main result by establishing Lemmas~\ref{Lem:Key1}, \ref{Lem:Key2}, and \ref{Lem:Key3}. We first generalize the construction of transition probabilities to the case of more general fillings of skyline diagrams, and provide an interpretation in terms of residues of certain rational functions. Although this level of generality is not strictly necessary for the proof of the main theorem, we formulate the definition in this setting as it appears to be the natural context for these objects. 

Our main technical tool for the proofs is Corollary~\ref{Cor:Varphi_Useful}, which describes the change resulting from modifying exactly one component of $\delta$.  We present the proof of Lemma~\ref{Lem:Key1} in Section~\ref{Sec:Key1}, Lemma~\ref{Lem:Key2} in Section~\ref{Sec:Key2}, and Lemma~\ref{Lem:Key3} in Section~\ref{Sec:Key3}.

\subsection{Preliminaries}

For technical reasons, we first generalize the previous definitions of $\boldsymbol{\delta}^{(r)}(T)$, $\Omega_r$, and $\mu_m(T)$ to more general fillings of skyline diagrams in $\bb{Z}_{>0}\times\bb{Z}_{>0}$. 

\begin{dfn}
A \emph{skyline diagram} $D\subset\bb{Z}_{>0}\times\bb{Z}_{>0}$ of size $n$ is a subset of the form
\begin{align*}
	D=\{(i,j)\in\bb{Z}_{>0}\times\bb{Z}_{>0}\mid 1\leq j \leq d_i\}
\end{align*}
for some sequence $(d_1,d_2,\ldots)$ of nonnegative integers such that $\sum_{i}d_i=n$. 
\end{dfn}

\begin{dfn}
A \emph{standard column-strict filling} of a skyline diagram $D$ is a bijection $T:D\rightarrow\{1,\ldots,|D|\}$ which is increasing along each column. We denote by $\msf{SSF}_n$ the set of standard column-strict fillings of skyline diagrams of size $n$, and by $\widetilde{\msf{V}}_n$ the $\bb{Q}(q)$-vector space formally spanned by $\msf{SSF}_n$.
\end{dfn}

\begin{dfn}
For $T\in\msf{SSF}_n$ and $r\in\bb{Z}_{\geq0}$, we define $\boldsymbol{\delta}^{(r)}(T)=(\delta_i)_{i\in\bb{Z}}\in\msf{Maya}$ as follows:
\begin{itemize}
	\item If $i\leq0$, then we set $\delta_i=1$.
	\item If $i>0$ and there are no boxes in the $i$-th column of $T$, then we set $\delta_i=0$.
	\item Otherwise, let $\msf{top}_T(i)$ be the entry in the topmost box of the $i$-th column of $T$. If $\msf{top}_T(i)>r$, then we set $\delta_i=1$. Otherwise, we set $\delta_i=0$.
\end{itemize}
\end{dfn}

\begin{dfn}
For $r\in\bb{Z}_{\geq0}$, we define a $\bb{Q}(q)$-linear map $\widetilde{\Omega}_r:\widetilde{\msf{V}}_{n-1}\rightarrow\widetilde{\msf{V}}_n$ by
\begin{align*}
	\widetilde{\Omega}_r(T)=\sum_{c\in\bb{Z}_{>0}}\varphi_c\bigl(\boldsymbol{\delta}^{(r)}(T);q\bigr) \mf{f}_c(T)
\end{align*}
for any $T\in\msf{SSF}_{n-1}$, where $\mf{f}_c(T)\in\msf{SSF}_n$ is obtained by adding $\boxed{n}$ on top of the $c$-th column of $T$. 
\end{dfn}

\begin{remark}
With respect to the obvious inclusion $\msf{V}_n\subset\widetilde{\msf{V}}_n$, $\widetilde{\Omega}_r$ restricts to the map $\Omega_r:\msf{V}_{n-1}\rightarrow\msf{V}_n$. One of the advantages of this generalization is that the condition $c\in W\bigl(\boldsymbol{\delta}^{(r)}(T)\bigr)$ is no longer required to ensure the well-definedness of $\mf{f}_c(T)$. 
\end{remark}

\begin{dfn}
For $T\in\msf{SSF}_n$ and $m\in[n-1]$, we define $\mu_m(T)\in\widetilde{\msf{V}}_n$ by 
\begin{align*}
	\mu_m(T)\coloneqq\begin{cases}
		\frac{[a-b-1]_q}{[a-b]_q}\cdot T+\frac{[b-a-1]_q}{[b-a]_q}\cdot\tau_m(T) &\text{if }\msf{top}_T(a)=m \text{ and }\msf{top}_T(b)=m+1,\\
		0 &\text{otherwise}.
	\end{cases}
\end{align*}
Here, $\tau_m(T)\in\msf{SSF}_n$ is the filling obtained by swapping the boxes containing $m$ and $m+1$.
\end{dfn}

We note that $\tau_m(T)$ is always well-defined as an element of $\msf{SSF}_n$ under the assumption that $\msf{top}_T(a)=m$ and $\msf{top}_T(b)=m+1$.

\subsection{Residues}

In this section, to avoid a lengthy case-by-case analysis in the proof of the Lemma~\ref{Lem:Key2}, we reinterpret the definition of $\varphi_c(\delta;q)$ in terms of residues of rational functions.

\begin{dfn}
Let $q\in\bb{R}_{>0}\setminus\{1\}$ and $\delta\in\msf{Maya}$ be a Maya diagram. We define a rational function $\psi_q(\delta;z)$ of $z$ by
\begin{align*}
	\psi_q(\delta;z)\coloneqq\frac{\prod_{b\in R(\delta)}(z-q^{-b})}{\prod_{a\in W(\delta)}(z-q^{-a})}.
\end{align*}
\end{dfn}

\begin{remark}\label{Rem:Res}
The function $\psi_q(\delta;z)$ has a simple pole at $z=q^{-c}$ for each $c\in W(\delta)$ and its residue is given by
\begin{align*}
	\Res_{z=q^{-c}}\psi_q(\delta;z)dz=\frac{\prod_{b\in R(\delta)}(q^{-c}-q^{-b})}{\prod_{a\in W(\delta)\setminus\{c\}}(q^{-c}-q^{-a})}=\varphi_c(\delta;q).
\end{align*}
We may consider $\psi_q(\delta;z)dz$ as a meromorphic differential on $\bb{P}^1$. Since the sum of all residues of a meromorphic differential on a compact Riemann surface is zero, we obtain
\begin{align*}
	\sum_{c\in W(\delta)}\varphi_c(\delta;q)=-\Res_{z=\infty}\psi_q(\delta;z)dz=1
\end{align*}
for any $q\in\bb{R}_{>0}\setminus\{1\}$. This implies $\sum_{c\in W(\delta)}\varphi_c(\delta;q)=1$ as a rational function in $q$.
\end{remark}

\begin{dfn}
Let $\delta=(\delta_i)_{i\in\bb{Z}}\in\msf{Maya}$ be a Maya diagram. For each $d\in\bb{Z}$ such that $\delta_d=0$ (resp. $\delta_d=1$), we denote by $\mf{f}_d(\delta)$ (resp. $\mf{e}_d(\delta)$) the Maya diagram obtained by replacing the $d$-th component of $\delta$ with $1$ (resp. $0$).
\end{dfn}

\begin{lemma}\label{Lem:Psi_Useful}
Let $\delta\in\msf{Maya}$ be a Maya diagram, $q\in\bb{R}_{>0}$, and $d\in\bb{Z}$ with $\delta_d=1$. We have
\begin{align*}
	\psi_q\left(\mf{e}_d(\delta);z\right)=\frac{z-q^{-d-1}}{z-q^{-d}}\psi_q(\delta;z).
\end{align*}
\end{lemma}

\begin{proof}
We divide the proof into four cases depending on the values of $\delta_{d-1}$ and $\delta_{d+1}$.

\medskip
\noindent \textbf{Case 1: $\delta_{d-1} = 0$ and $\delta_{d+1} = 0$.}

In this case, we have $\delta=(\ldots,\overset{d-1}{0},\overset{d}{1},\overset{d+1}{0},\ldots)$ and $\mf{e}_d(\delta)=(\ldots,\overset{d-1}{0},\overset{d}{0},\overset{d+1}{0},\ldots)$. This implies
\begin{align*}
	W\left(\mf{e}_d(\delta)\right)&=W(\delta)\setminus\{d+1\},\\
	R\left(\mf{e}_d(\delta)\right)&=R(\delta)\setminus\{d\}.
\end{align*}
Hence we obtain
\begin{align*}
	\psi_q(\mf{e}_d(\delta);z)=\frac{\prod_{b\in R(\delta)\setminus\{d\}}(z-q^{-b})}{\prod_{a\in W(\delta)\setminus\{d+1\}}(z-q^{-a})}=\frac{z-q^{-d-1}}{z-q^{-d}}\psi_q(\delta;z).
\end{align*}

\medskip
\noindent \textbf{Case 2: $\delta_{d-1} = 1$ and $\delta_{d+1} = 0$.}

In this case, we have $\delta=(\ldots,\overset{d-1}{1},\overset{d}{1},\overset{d+1}{0},\ldots)$ and $\mf{e}_d(\delta)=(\ldots,\overset{d-1}{1},\overset{d}{0},\overset{d+1}{0},\ldots)$. This implies
\begin{align*}
	W\left(\mf{e}_d(\delta)\right)&=W(\delta)\cup\{d\}\setminus\{d+1\},\\
	R\left(\mf{e}_d(\delta)\right)&=R(\delta).
\end{align*}
Hence we obtain
\begin{align*}
	\psi_q(\mf{e}_d(\delta);z)=\frac{\prod_{b\in R(\delta)}(z-q^{-b})}{\prod_{a\in W(\delta)\cup\{d\}\setminus\{d+1\}}(z-q^{-a})}=\frac{z-q^{-d-1}}{z-q^{-d}}\psi_q(\delta;z).
\end{align*}

\medskip
\noindent \textbf{Case 3: $\delta_{d-1} = 0$ and $\delta_{d+1} = 1$.}

In this case, we have $\delta=(\ldots,\overset{d-1}{0},\overset{d}{1},\overset{d+1}{1},\ldots)$ and $\mf{e}_d(\delta)=(\ldots,\overset{d-1}{0},\overset{d}{0},\overset{d+1}{1},\ldots)$. This implies
\begin{align*}
	W\left(\mf{e}_d(\delta)\right)&=W(\delta),\\
	R\left(\mf{e}_d(\delta)\right)&=R(\delta)\cup\{d+1\}\setminus\{d\}.
\end{align*}
Hence we obtain
\begin{align*}
	\psi_q(\mf{e}_d(\delta);z)=\frac{\prod_{b\in R(\delta)\cup\{d+1\}\setminus\{d\}}(z-q^{-b})}{\prod_{a\in W(\delta)}(z-q^{-a})}=\frac{z-q^{-d-1}}{z-q^{-d}}\psi_q(\delta;z).
\end{align*}

\medskip
\noindent \textbf{Case 4: $\delta_{d-1} = 1$ and $\delta_{d+1} = 1$.}

In this case, we have $\delta=(\ldots,\overset{d-1}{1},\overset{d}{1},\overset{d+1}{1},\ldots)$ and $\mf{e}_d(\delta)=(\ldots,\overset{d-1}{1},\overset{d}{0},\overset{d+1}{1},\ldots)$. This implies
\begin{align*}
	W\left(\mf{e}_d(\delta)\right)&=W(\delta)\cup\{d\},\\
	R\left(\mf{e}_d(\delta)\right)&=R(\delta)\cup\{d+1\}.
\end{align*}
Hence we obtain
\begin{align*}
	\psi_q(\mf{e}_d(\delta);z)=\frac{\prod_{b\in R(\delta)\cup\{d+1\}}(z-q^{-b})}{\prod_{a\in W(\delta)\cup\{d\}}(z-q^{-a})}=\frac{z-q^{-d-1}}{z-q^{-d}}\psi_q(\delta;z).
\end{align*}
\end{proof}

\begin{corollary}\label{Cor:Varphi_Useful}
Let $\delta\in\msf{Maya}$ be a Maya diagram, $d\in\bb{Z}$ with $\delta_d=1$, and $c\neq d$. We have
\begin{align*}
	\varphi_c(\mf{e}_d(\delta);q)=\frac{[c-d-1]_q}{[c-d]_q}\varphi_c(\delta;q).
\end{align*}
\end{corollary}

\begin{proof}
Since the assertion is an identity of rational functions in $q$, it suffices to check it for any $q\in\bb{R}_{>0}\setminus\{1\}$. By Remark~\ref{Rem:Res}, Lemma~\ref{Lem:Psi_Useful}, and $c\neq d$, we obtain
\begin{align*}
	\varphi_c(\mf{e}_d(\delta);q)&=\Res_{z=q^{-c}}\left(\frac{z-q^{-d-1}}{z-q^{-d}}\psi_{q}(\delta;z)dz\right)\\
	&=\frac{q^{-c}-q^{-d-1}}{q^{-c}-q^{-d}}\Res_{z=q^{-c}}\left(\psi_{q}(\delta;z)dz\right)\\
	&=\frac{[c-d-1]_q}{[c-d]_q}\varphi_c(\delta;q).
\end{align*}
\end{proof}

\subsection{Proof of Lemma~\ref{Lem:Key1}}\label{Sec:Key1}

In this section, we verify Lemma~\ref{Lem:Key1}. Since the restriction of $\widetilde{\Omega}_r$ to the subspace $\msf{V}_{n}\subset\widetilde{\msf{V}}_{n}$ coincides with $\Omega_r$, Lemma~\ref{Lem:Key1} follows from the following more general statement.

\begin{prop}\label{Prop:Key1}
For any $T\in\msf{SSF}_{n-1}$ and $0\leq r<n$, we have
\begin{align*}
	\widetilde{\Omega}_r\widetilde{\Omega}_r(T)=\sum_{\substack{c< d\in\bb{Z}_{>0}\\\delta_c=\delta_d=0}}\varphi_c\bigl(\mf{f}_d(\boldsymbol{\delta}^{(r)}(T));q\bigr)\varphi_d\bigl(\mf{f}_c(\boldsymbol{\delta}^{(r)}(T));q\bigr)\mu_n\bigl(\mf{f}_d\mf{f}_c(T)\bigr).
\end{align*}
\end{prop}

\begin{proof}
Let $\delta=\boldsymbol{\delta}^{(r)}(T)$. For any $c\in \bb{Z}_{>0}$ with $\delta_c=0$, we have $\boldsymbol{\delta}^{(r)}(\mf{f}_c(T))=\mf{f}_c(\delta)$, since the added entry $n$ satisfies $n>r$. Since $\varphi_c(\delta)=0$ unless $\delta_c=0$, and $\varphi_c(\mf{f}_c(\delta);q)=0$, we obtain
\begin{align*}
	\widetilde{\Omega}_r\widetilde{\Omega}_r(T)&=\widetilde{\Omega}_r\left(\sum_{\substack{c\in\bb{Z}_{>0}\\\delta_c=0}}\varphi_c(\delta;q)\mf{f}_c(T)\right)\\
&=\sum_{\substack{c\neq d\in\bb{Z}_{>0}\\\delta_c=\delta_d=0}}\varphi_c(\delta;q)\varphi_d(\mf{f}_c(\delta);q)\mf{f}_d\mf{f}_c(T)\\
&=\sum_{\substack{c\neq d\in\bb{Z}_{>0}\\\delta_c=\delta_d=0}}\varphi_c(\mf{f}_d(\delta);q)\varphi_d(\mf{f}_c(\delta);q)\frac{[c-d-1]_q}{[c-d]_q}\mf{f}_d\mf{f}_c(T)\\
&=\sum_{\substack{c< d\in\bb{Z}_{>0}\\\delta_c=\delta_d=0}}\varphi_c(\mf{f}_d(\delta);q)\varphi_d(\mf{f}_c(\delta);q)\left(\frac{[c-d-1]_q}{[c-d]_q}\mf{f}_d\mf{f}_c(T)+\frac{[d-c-1]_q}{[d-c]_q}\mf{f}_c\mf{f}_d(T)\right)\\
&=\sum_{\substack{c< d\in\bb{Z}_{>0}\\\delta_c=\delta_d=0}}\varphi_c(\mf{f}_d(\delta);q)\varphi_d(\mf{f}_c(\delta);q)\mu_n\bigl(\mf{f}_d\mf{f}_c(T)\bigr),
\end{align*}
where we used Corollary~\ref{Cor:Varphi_Useful} for $\mf{f}_d(\delta)$ and $d\in\bb{Z}$ in the third equality.
\end{proof}

\begin{proof}[Proof of Lemma~\ref{Lem:Key1}]

Proposition~\ref{Prop:Key1} implies that for any $T\in\msf{SYT}_{n-1}$, with $\delta=\boldsymbol{\delta}^{(r)}(T)$, we have
\begin{align*}
	\Omega_r\Omega_r(T)=\sum_{\substack{c< d\in\bb{Z}_{>0}\\\delta_c=\delta_d=0}}\varphi_c(\mf{f}_d(\delta);q)\varphi_d(\mf{f}_c(\delta);q)\mu_n\bigl(\mf{f}_d\mf{f}_c(T)\bigr).
\end{align*}
We note that $c<d$ implies $\frac{[c-d-1]_q}{[c-d]_q}\neq0$ and hence $\mf{f}_d\mf{f}_c(T)$ appears in $\mu_n(\mf{f}_d\mf{f}_c(T))$ with a non-zero coefficient. Since the fillings of the form $\mf{f}_d\mf{f}_c(T)$ and $\mf{f}_c\mf{f}_d(T)$ are pairwise distinct as $c$ and $d$ vary over positive integers with $c<d$, and since $\Omega_r\Omega_r(T)$ must lie in $\msf{V}_{n+1}$, the coefficient $\varphi_c(\mf{f}_d(\delta);q)\varphi_d(\mf{f}_c(\delta);q)$ of $\mu_n(\mf{f}_d\mf{f}_c(T))$ must be zero unless $\mf{f}_d\mf{f}_c(T)\in\msf{SYT}_{n+1}$. Consequently, every non-zero term in the sum is of the form $\mu_n(T')$ for some $T'\in\msf{SYT}_{n+1}$. This proves $\Omega_r\Omega_r(T)\in\msf{M}_{n,n+1}$. 
\end{proof}

\subsection{Proof of Lemma~\ref{Lem:Key2}}\label{Sec:Key2}

In this section, we verify Lemma~\ref{Lem:Key2}. It suffices to prove the following more general result.

\begin{prop}\label{Prop:Key2}
For any $T\in\msf{SSF}_n$ and $r\in[n-1]$, we have
\begin{align*}
\left(\widetilde{\Omega}_r-\frac{\widetilde{\Omega}_{r-1}+q\widetilde{\Omega}_{r+1}}{1+q}\right)\bigl(\mu_r(T)\bigr)\in\Span_{\bb{Q}(q)}\left\{\mf{f}_c(T)-\tau_r\bigl(\mf{f}_c(T)\bigr)\relmiddle| c\in\bb{Z}_{>0}\right\}
\end{align*}
\end{prop}

\begin{proof}
Let $T\in\msf{SSF}_n$ satisfy $\msf{top}_T(a)=r$ and $\msf{top}_T(b)=r+1$. We set $\delta=\boldsymbol{\delta}^{(r)}(T)=(\ldots,\overset{a}{0},\ldots,\overset{b}{1},\ldots)$. We then have 
\begin{gather*}
	\boldsymbol{\delta}^{(r-1)}(T)=\mf{f}_a(\delta)=\boldsymbol{\delta}^{(r-1)}(\tau_r(T)),\\
	\boldsymbol{\delta}^{(r+1)}(T)=\mf{e}_b(\delta)=\boldsymbol{\delta}^{(r+1)}(\tau_r(T)),\\
	\boldsymbol{\delta}^{(r)}(\tau_r(T))=\mf{f}_a\mf{e}_b(\delta).
\end{gather*}
Hence we obtain
\begin{align*}
&\left(\widetilde{\Omega}_r-\frac{\widetilde{\Omega}_{r-1}+q\widetilde{\Omega}_{r+1}}{1+q}\right)\bigl(\mu_r(T)\bigr)\\
&\qquad=\sum_{c\in\bb{Z}_{>0}}\frac{[a-b-1]_q}{[a-b]_q}\left(\varphi_c(\delta;q)-\frac{\varphi_c\left(\mf{f}_a(\delta);q\right)+q\varphi_{c}(\mf{e}_b(\delta);q)}{1+q}\right)\mf{f}_c(T)\\
&\qquad\quad+\sum_{c\in\bb{Z}_{>0}}\frac{[b-a-1]_q}{[b-a]_q}\left(\varphi_c(\mf{f}_a\mf{e}_b(\delta);q)-\frac{\varphi_c\left(\mf{f}_a(\delta);q\right)+q\varphi_{c}(\mf{e}_b(\delta);q)}{1+q}\right)\mf{f}_c(\tau_r(T)).
\end{align*}
 
It suffices to check that the sum of the coefficient of $\mf{f}_c(T)$ and the coefficient of $\mf{f}_c(\tau_r(T))$ is zero for each $c\in\bb{Z}$. Since this is an equality of rational functions in $q$, it is enough to check it for $q\in\bb{R}_{>0}\setminus\{1\}$. 

By Remark~\ref{Rem:Res} and Lemma~\ref{Lem:Psi_Useful}, we may express the coefficient of $\mf{f}_c(T)$ as the residue at $z=q^{-c}$ of
\begin{align*}
&\frac{[a-b-1]_q}{[a-b]_q}\left(\frac{z-q^{-a-1}}{z-q^{-a}}-\frac{1}{1+q}-\frac{q}{1+q}\frac{z-q^{-b-1}}{z-q^{-b}}\frac{z-q^{-a-1}}{z-q^{-a}}\right)\psi_q(\mf{f}_a(\delta);z),\\
&\qquad=\frac{(1-q)z(q^{-a}-q^{-b-1})(q^{-b}-q^{-a-1})}{(1+q)(z-q^{-a})(z-q^{-b})(q^{-a}-q^{-b})}\psi_q(\mf{f}_a(\delta);z).
\end{align*}

On the other hand, the coefficient of $\mf{f}_c(\tau_r(T))$ is the residue at $z=q^{-c}$ of the rational function obtained by swapping $a$ and $b$ in the coefficient multiplying $\psi_q(\mf{f}_a(\delta);z)$ above, i.e.,
\begin{align*}
&\frac{[b-a-1]_q}{[b-a]_q}\left(\frac{z-q^{-b-1}}{z-q^{-b}}-\frac{1}{1+q}-\frac{q}{1+q}\frac{z-q^{-b-1}}{z-q^{-b}}\frac{z-q^{-a-1}}{z-q^{-a}}\right)\psi_q(\mf{f}_a(\delta);z)\\
&\qquad=-\frac{(1-q)z(q^{-a}-q^{-b-1})(q^{-b}-q^{-a-1})}{(1+q)(z-q^{-a})(z-q^{-b})(q^{-a}-q^{-b})}\psi_q(\mf{f}_a(\delta);z).
\end{align*}
Therefore, the sum of these rational functions vanishes identically, implying that the sum of their residues is zero.
\end{proof}

\begin{proof}[Proof of Lemma~\ref{Lem:Key2}]

Let $T\in\msf{SYT}_n$ be a standard Young tableau such that $\mu_r(T)\neq0$. Proposition~\ref{Prop:Key2} implies that we may write
\begin{align*}
	\left(\Omega_r-\frac{\Omega_{r-1}+q\Omega_{r+1}}{1+q}\right)\bigl(\mu_r(T)\bigr)=\sum_{c\in\bb{Z}_{>0}}g_c(q)\bigl(\mf{f}_c(T)-\tau_r(\mf{f}_c(T))\bigr)
\end{align*}
for some $g_c(q)\in\bb{Q}(q)$. Since the LHS belongs to $\msf{V}_{n+1}$, and the fillings of the form $\mf{f}_c(T)$ and $\tau_r(\mf{f}_c(T))$ are pairwise distinct as $c$ varies over positive integers, the coefficient $g_c(q)$ must vanish unless $\mf{f}_c(T)$ and $\tau_r(\mf{f}_c(T))$ are standard Young tableaux. Thus, the RHS lies in $\msf{K}_{r,n+1}$.
\end{proof}

\subsection{Proof of Lemma~\ref{Lem:Key3}}\label{Sec:Key3}

In this section, we prove Lemma~\ref{Lem:Key3} and hence complete the proof of our main theorem. 

\begin{prop}\label{Prop:Key3}
For any $T\in\msf{SSF}_{n-1}$ and $0\leq r<n-1$, we have
\begin{align*}
\left(\widetilde{\Omega}_{r+1}\widetilde{\Omega}_r-\frac{\widetilde{\Omega}_{r}\widetilde{\Omega}_r+q\widetilde{\Omega}_{r+1}\widetilde{\Omega}_{r+1}}{1+q}\right)(T)\in\Span_{\bb{Q}(q)}\left\{\mf{f}_c\mf{f}_d(T)-\mf{f}_d\mf{f}_c(T)\relmiddle|c<d\in\bb{Z}_{>0}\right\}.
\end{align*}
\end{prop}

\begin{proof}

Let $\delta=\boldsymbol{\delta}^{(r)}(T)$. If the box labeled by $r+1$ does not appear at the top of some column in $T$, then we have $\boldsymbol{\delta}^{(r+1)}(T)=\delta$ and $\boldsymbol{\delta}^{(r)}(\mf{f}_c(T))=\boldsymbol{\delta}^{(r+1)}(\mf{f}_c(T))$ for any $c\in\bb{Z}_{>0}$. This implies 
\begin{align*}
	\widetilde{\Omega}_{r}\widetilde{\Omega}_r(T)=\widetilde{\Omega}_{r+1}\widetilde{\Omega}_r(T)=\widetilde{\Omega}_{r+1}\widetilde{\Omega}_{r+1}(T).
\end{align*}
Thus, the assertion holds trivially in this case.

Therefore, we may assume that $\msf{top}_T(b)=r+1$ for some $b\in\bb{Z}_{>0}$. In this case, we have $\delta_b=1$ and $\boldsymbol{\delta}^{(r+1)}(T)=\mf{e}_b(\delta)$. Moreover, if $c\in\bb{Z}_{>0}$ satisfies $\delta_c=0$, then we have $\boldsymbol{\delta}^{(r+1)}(\mf{f}_c(T))=\mf{f}_c\mf{e}_b(\delta)$ and its $d$-th component is zero if $d=b$ or $\delta_d=0$ with $d\neq c$. Hence we obtain
\begin{align*}
	\widetilde{\Omega}_{r+1}\widetilde{\Omega}_r(T)&=\widetilde{\Omega}_{r+1}\left(\sum_{\substack{c\in\bb{Z}_{>0}\\\delta_c=0}}\varphi_c(\delta;q)\mf{f}_c(T)\right)\\
	&=\sum_{\substack{c\in\bb{Z}_{>0}\\\delta_c=0}}\varphi_c(\delta;q)\varphi_b(\mf{f}_c\mf{e}_b(\delta);q)\mf{f}_b\mf{f}_c(T)+\sum_{\substack{c\neq d\in\bb{Z}_{>0}\\\delta_c=\delta_d=0}}\varphi_c(\delta;q)\varphi_d(\mf{f}_c\mf{e}_b(\delta);q)\mf{f}_d\mf{f}_c(T)
\end{align*}
On the other hand, we have
\begin{align*}
\widetilde{\Omega}_{r}\widetilde{\Omega}_r(T)&=\sum_{\substack{c\neq d\in\bb{Z}_{>0}\\\delta_c=\delta_d=0}}\varphi_c(\delta;q)\varphi_d(\mf{f}_c(\delta);q)\mf{f}_d\mf{f}_c(T),\\
\widetilde{\Omega}_{r+1}\widetilde{\Omega}_{r+1}(T)&=\sum_{\substack{c\in\bb{Z}_{>0}\\\delta_c=0}}\varphi_c(\mf{e}_b(\delta);q)\varphi_b(\mf{f}_c\mf{e}_b(\delta);q)\mf{f}_b\mf{f}_c(T)+\sum_{\substack{c\in\bb{Z}_{>0}\\\delta_c=0}}\varphi_b(\mf{e}_b(\delta);q)\varphi_c(\delta;q)\mf{f}_c\mf{f}_b(T)\\
&\quad+\sum_{\substack{c\neq d\in\bb{Z}_{>0}\\\delta_c=\delta_d=0}}\varphi_c(\mf{e}_b(\delta);q)\varphi_d(\mf{f}_c\mf{e}_b(\delta);q)\mf{f}_d\mf{f}_c(T).
\end{align*}
Therefore, by repeatedly applying Corollary~\ref{Cor:Varphi_Useful}, 
we obtain
\begin{align*}
&\left(\widetilde{\Omega}_{r+1}\widetilde{\Omega}_r-\frac{\widetilde{\Omega}_{r}\widetilde{\Omega}_r+q\widetilde{\Omega}_{r+1}\widetilde{\Omega}_{r+1}}{1+q}\right)(T)\\
&\qquad=\sum_{\substack{c\in\bb{Z}_{>0}\\\delta_c=0}}\left(\varphi_c(\delta;q)\varphi_b(\mf{f}_c\mf{e}_b(\delta);q)-\frac{q}{1+q}\varphi_c(\mf{e}_b(\delta);q)\varphi_{b}(\mf{f}_c\mf{e}_b(\delta);q)\right)\mf{f}_b\mf{f}_c(T)\\
&\qquad\quad-\sum_{\substack{c\in\bb{Z}_{>0}\\\delta_c=0}}\frac{q}{1+q}\varphi_b(\mf{e}_b(\delta);q)\varphi_c(\delta;q)\mf{f}_c\mf{f}_b(T)\\
&\qquad\quad+\sum_{\substack{c\neq d\in\bb{Z}_{>0}\\\delta_c=\delta_d=0}}\left(\varphi_c(\delta;q)\varphi_d(\mf{f}_c\mf{e}_b(\delta);q)-\frac{\varphi_c(\delta;q)\varphi_d(\mf{f}_c(\delta);q)+q\varphi_c(\mf{e}_b(\delta);q)\varphi_d(\mf{f}_c\mf{e}_b(\delta);q)}{1+q}\right)\mf{f}_d\mf{f}_c(T).\\
\end{align*}

It suffices to check that the coefficients of $\mf{f}_d\mf{f}_c(T)$ is antisymmetric with respect to $c$ and $d$ for each pair $c\neq d\in\{i\in\bb{Z}\mid\delta_i=0\}\cup\{b\}$. By Corollary~\ref{Cor:Varphi_Useful}, the coefficient of $\mf{f}_b\mf{f}_c(T)$ is given by
\begin{align*}
	\left(1-\frac{q}{1+q}\frac{[c-b-1]_q}{[c-b]_q}\right)\frac{[b-c]_q}{[b-c-1]_q}\varphi_c(\delta;q)\varphi_b(\mf{e}_b(\delta);q)=\frac{q}{1+q}\varphi_c(\delta;q)\varphi_b(\mf{e}_b(\delta);q),
\end{align*}
which is the negative of the coefficient of $\mf{f}_c\mf{f}_b(T)$.

The coefficient of $\mf{f}_d\mf{f}_c(T)$ for $c\neq d\in\{i\in\bb{Z}_{>0}\mid\delta_i=0\}$ is given by
\begin{align*}
&\left(\frac{[d-b-1]_q}{[d-b]_q}-\frac{1}{1+q}-\frac{q}{1+q}\frac{[c-b-1]_q}{[c-b]_q}\frac{[d-b-1]_q}{[d-b]_q}\right)\frac{[c-d-1]_q}{[c-d]_q}\varphi_c(\mf{f}_d(\delta);q)\varphi_d(\mf{f}_c(\delta);q)\\
&\qquad=\frac{(1-q)q^b(q^c-q^{d-1})(q^{c-1}-q^{d})}{(1+q)(q^b-q^c)(q^b-q^d)(q^c-q^d)}\varphi_c(\mf{f}_d(\delta);q)\varphi_d(\mf{f}_c(\delta);q).
\end{align*}
Since the last expression is antisymmetric in $c$ and $d$, the assertion follows.
\end{proof}

\begin{proof}[Proof of Lemma~\ref{Lem:Key3}]
Let $T\in\msf{SYT}_{n-1}$ and $0\leq r<n-1$. Proposition~\ref{Prop:Key3} implies that 
\begin{align*}
\left(\Omega_{r+1}\Omega_r-\frac{\Omega_r\Omega_r+q\Omega_{r+1}\Omega_{r+1}}{1+q}\right)(T)=\sum_{c<d\in\bb{Z}_{>0}}g_{c,d}(q)\bigl(\mf{f}_d\mf{f}_c(T)-\mf{f}_c\mf{f}_d(T)\bigr)\in\msf{V}_{n+1}
\end{align*}
for some $g_{c,d}(q)\in\bb{Q}(q)$. Since the fillings of the form $\mf{f}_c\mf{f}_d(T)$ are pairwise distinct as $c$ and $d$ vary over positive integers with $c\neq d$, the coefficient $g_{c,d}(q)$ must be zero unless $\mf{f}_d\mf{f}_c(T),\mf{f}_c\mf{f}_d(T)\in\msf{SYT}_{n+1}$. Therefore, the RHS lies in $\msf{K}_{n,n+1}$.
\end{proof}

\appendix

\section{Review of probability theory}\label{Sec:Prob}

As some readers of this paper may not be familiar with probability theory, we briefly recall the basic language of probability theory in this appendix. For more details, see for example \cite{Cin11}.

\subsection{Stochastic process}

The fundamental setting in probability theory is a probability space $(\Omega, \msc{H}, \Prob)$. Here, $\Omega$ is a set, $\msc{H}$ is a $\sigma$-algebra on $\Omega$ (i.e., a non-empty collection of subsets of $\Omega$ closed under complements and countable unions), and $\Prob$ is a measure on the measurable space $(\Omega,\msc{H})$ such that $\Prob(\Omega)=1$. Such a probability space is regarded as a mathematical model of a random experiment. 

For simplicity, we only consider stochastic processes with a countable state space and parameter set $\bb{Z}_{\geq0}$ (or $\{0,1,\ldots,n\}$) here. Let $S$ be a countable set equipped with the $\sigma$-algebra $\msc{S}$ consisting of all subsets of $S$. A map $X:\Omega\rightarrow S$ is called a \emph{random variable} if it is measurable, i.e., $X^{-1}(A)\in\msc{H}$ for any $A\in\msc{S}$. This induces a probability measure on $(S,\msc{S})$, called the \emph{distribution} of $X$, given by sending $A\in\msc{S}$ to the probability that $X$ is in $A$, i.e., 
\begin{align*}
	A\mapsto\Prob\{X\in A\}\coloneqq\Prob(X^{-1}(A))
\end{align*}
for any $A\in\msc{S}$. We note that a measurable map $f:S\rightarrow S'$ defines a random variable $f\circ X:\Omega\rightarrow S'$ via composition.

A collection $\{X_n\}_{n\in \bb{Z}_{\geq0}}$ of random variables $X_n:\Omega\rightarrow S$ indexed by $n\in \bb{Z}_{\geq0}$ is called a \emph{stochastic process} with \emph{state space} $(S,\msc{S})$ and \emph{parameter set} $\bb{Z}_{\geq0}$. The \emph{probability law} of the stochastic process $\{X_n\}_{n\in \bb{Z}_{\geq0}}$ is the data which assigns the value
\begin{align*}
	\Prob\{X_{i_1}\in A_1,\ldots,X_{i_n}\in A_n\}\coloneqq\Prob\left(X_{i_1}^{-1}(A_1)\cap\cdots\cap X_{i_n}^{-1}(A_n)\right)
\end{align*}
for any $n\in\bb{Z}_{>0}$, $i_1,\ldots,i_n\in \bb{Z}_{\geq0}$, and $A_1,\ldots,A_n\in\msc{S}$. This data is usually enough to calculate probabilistic quantities related to the stochastic process $\{X_n\}_{n\in \bb{Z}_{\geq0}}$.

\subsection{Probability transition kernels}

In order to construct a stochastic process, it is useful to consider its probability transition kernels. For each $n\in\bb{Z}_{\geq0}$, we consider the \emph{probability transition kernel} $K_{n}:S^n\times\mca{S}\rightarrow[0,1]$ given by
\begin{align*}
	K_n(s_0,\ldots,s_{n-1};B)&=\Prob\{X_n\in B\mid X_0=s_0,\ldots,X_{n-1}=s_{n-1}\}\\
	&\coloneqq\frac{\Prob\{X_0=s_0,\ldots,X_{n-1}=s_{n-1},X_n\in B\}}{\Prob\{X_0=s_0,\ldots,X_{n-1}=s_{n-1}\}}
\end{align*}
for any $(s_0,\ldots,s_{n-1})\in S^n$ such that $\Prob\{X_0=s_0,\ldots,X_{n-1}=s_{n-1}\}\neq0$. When $n=0$, we regard $K_0:\msc{S}\rightarrow[0,1]$ as the distribution of $X_0$. If $\Prob\{X_0=s_0,\ldots,X_{n-1}=s_{n-1}\}=0$, then we may consider any probability measure $\mu:\msc{S}\rightarrow[0,1]$ and set $K_n(s_0,\ldots,s_{n-1};B)=\mu(B)$. These data determine the probability law of $\{X_n\}_{n\in \bb{Z}_{\geq0}}$ since we have, for example,
\begin{align*}
	\Prob\{X_0=s_0,\ldots,X_n=s_n\}=K_0\left(\{s_0\}\right)K_1\left(s_0;\{s_1\}\right)K_2\left(s_0,s_1;\{s_2\}\right)\cdots K_n\left(s_0,\ldots,s_{n-1};\{s_n\}\right).
\end{align*}

Conversely, any collection of maps $K_{n}:S^n\times\mca{S}\rightarrow[0,1]$ such that the map $B\mapsto K_n(s_0,\ldots,s_{n-1};B)$ is a probability measure for any fixed $s_0,\ldots,s_{n-1}\in S$ (and the map $(s_0,\ldots,s_{n-1})\mapsto K_n(s_0,\ldots,s_{n-1};B)$ is measurable for any fixed $B\in\msc{S}$ for a general measurable space $(S,\msc{S})$) determines a stochastic process $\{X_n\}_{n\in \bb{Z}_{\geq0}}$ on some probability space $(\Omega,\msc{H},\Prob)$ whose probability transition kernels coincide with $\{K_n\}_{n\in\bb{Z}_{\geq0}}$ (see \cite[Chapter~4]{Cin11} for more details). This allows us to construct a stochastic process by only specifying $\{K_n\}_{n\in\bb{Z}_{\geq0}}$, without explicitly constructing the probability space $(\Omega,\msc{H},\Prob)$.

\subsection{Markov chains}

A stochastic process $\{X_n\}_{n\in \bb{Z}_{\geq0}}$ is called a \emph{Markov chain} if the future behavior does not depend on the past history but only on the current state, i.e., there exists a \emph{Markov kernel} $P_{n}:S\times\msc{S}\rightarrow[0,1]$ for any $n\in\bb{Z}_{>0}$ such that
\begin{align*}
	K_n(s_0,\ldots,s_{n-1};B)=P_n(s_{n-1};B)
\end{align*}
for any $s_0,\ldots,s_{n-1}\in S$ and $B\in\msc{S}$. A Markov chain is called \emph{time-homogeneous} if its Markov kernel $P_n$ does not depend on $n$. Otherwise, it is called \emph{time-inhomogeneous}. For $s,s'\in S$, the value $P_n(s;\{s'\})$ is called the \emph{transition probability} from $s$ to $s'$ at step $n$.

\bibliography{qt-CSF}
\bibliographystyle{plain}

\end{document}